\documentclass[10pt]{amsart}
\usepackage{amscd,amssymb,amsmath,latexsym,enumerate}
\usepackage[mathscr]{euscript}
\usepackage{epsfig}
\newtheorem{theo}{Theorem}
\newtheorem{defini}{Definition}
\newtheorem{proposi}{Proposition}
\newtheorem{lemma}{Lemma}

\newtheorem{rem}{Remark}
\newtheorem{exam}{Example}

\newcommand{\Aa}{{\mathcal A}}
\newcommand{\Bb}{{\mathcal B}}
\newcommand{\Cc}{{\mathcal C}}
\newcommand{\Dd}{{\mathcal D}}
\newcommand{\Ee}{{\mathcal E}}

\newcommand{\Gg}{{\mathcal G}}
\newcommand{\Hh}{{\mathcal H}}

\newcommand{\Ll}{{\mathcal L}}

\newcommand{\Nn}{{\mathcal N}}
\newcommand{\Oo}{{\mathcal O}}

\newcommand{\Tt}{{\mathcal T}}

\newcommand{\Vv}{{\mathcal V}}
\newcommand{\Ww}{{\mathcal W}}

\newcommand{\RR}{{\bf R}}

\newcommand{\id}{{\mathbf 1}}

\newcommand{\CM}{{\mathbb C}}
\newcommand{\EM}{{\mathbb E}}

\newcommand{\NM}{{\mathbb N}}
\newcommand{\QM}{{\mathbb Q}}

\newcommand{\RM}{{\mathbb R}}

\newcommand{\ZM}{{\mathbb Z}}

\newcommand{\es}{{\mathscr E}}

\newcommand{\vs}{{\mathscr V}}

\newcommand{\TR}{{\rm Tr\,}}                       
\newcommand{\Cs}{$C^{\ast}$-algebra }              
\newcommand{\Css}{$C^{\ast}$-algebras }            
\newcommand{\pT}{\partial \Tt}                     
\newcommand{\diam}{\mbox{\rm diam}}                
\newcommand{\lip}{\Cc_{\mbox{\tiny\rm Lip}}}       
\newcommand{\theight}{{\mbox{\rm ht}}}             
\newcommand{\mheight}{{\mbox{\rm \tiny ht}}}       
\newcommand{\tdiam}{{\mbox{\rm diam}}}             
\newcommand{\tdist}{{\mbox{\rm dist}}}             
\newcommand{\ubox}{\overline{\mbox{\rm dim}}_B}    
\newcommand{\Dom}{\text {Dom}}      
\newcommand{\Supp}{\text {Supp}}     

\begin{document}

\title{Noncommutative Riemannian Geometry and Diffusion on Ultrametric Cantor Sets}
\thanks{Work supported in part by NSF grants DMS 06009565.  }

\author{John Pearson, Jean Bellissard}

\address{Georgia Institute of Technology\\  School of Mathematics\\ Atlanta, GA 30332-0160}
\email{pearson@math.gatech.edu, jeanbel@math.gatech.edu}

\begin{abstract}
An analogue of the Riemannian Geometry for an ultrametric Cantor set $(C,d)$ is described using the tools of Noncommutative Geometry. Associated with $(C,d)$ is a weighted rooted tree, its Michon tree \cite{MIC}. This tree allows to define a family of spectral triples $(\lip(C),\Hh,D)$ using the $\ell^2$-space of its vertices, giving the Cantor set the structure of a noncommutative Riemannian manifold. Here $\lip(C)$ denotes the space of Lipschitz continuous functions on $(C,d)$. The family of spectral triples is indexed by the space of {\em choice functions} which is shown to be the analogue of the sphere bundle of a Riemannian manifold. The Connes metric coming from the Dirac operator $D$ then allows to recover the metric on $C$. The corresponding $\zeta$-function is shown to have abscissa of convergence, $s_0$, equal to the {\em upper box dimension} of $(C,d)$. Taking the residue at this singularity leads to the definition of a canonical probability measure on $C$ which in certain cases coincides with the Hausdorff measure at dimension $s_0$. This measure in turns induces a measure on the space of choices. Given a choice, the commutator of $D$ with a Lipschitz continuous function can be intepreted as a directional derivative. By integrating over all choices, this leads to the definition of an analogue of the Laplace-Beltrami operator. This operator has compact resolvent and generates a Markov semigroup which plays the role of a Brownian motion on $C$. This construction is applied to the simplest case, the triadic Cantor set where: (i) the spectrum and the eigenfunctions of the Laplace-Beltrami operator are computed, (ii) the Weyl asymptotic formula is shown to hold with the dimension $s_0$, (iii) the corresponding Markov process is shown to have an anomalous diffusion with $\EM(d(X_{t},X_{t+\delta t})^2)\simeq \delta t\ln{(1/\delta t)}$ as $\delta t \downarrow 0$. 
\end{abstract}


\maketitle


\pagestyle{myheadings}
\markboth{John Pearson, Jean Bellissard}{Noncommutative Riemannian Geometry and Diffusion on Ultrametric Cantor Sets}


\tableofcontents



\vspace{1.5cm}

\section{Introduction}
\label{cantor07.sect-intro}

\noindent The present work aims to define a Noncommutative Riemannian structure on an ultrametric Cantor set $(C,d)$. To this end, according to the theory developed by Connes \cite{CON}, a spectral triple $(\Aa, \Hh,D)$ will be defined, namely $\Aa \subset \Cc(C)$ is a dense $\ast$-subalgebra of the \Css of complex valued continuous functions over $C$ , $\Hh$ is a Hilbert space on which the algebra $\Cc(C)$ is represented and $D$, called the Dirac operator, is a selfadjoint operator on $\Hh$ with compact resolvent, such that $[D,a]\in \Bb(\Hh)$ whenever $a\in\Aa$. In the present work, the spectral triple is {\em even}. Namely there is an operator $\Gamma$ on $\Hh$ such that $\Gamma=\Gamma^\ast$ and $\Gamma^2=\id$ and that $\Gamma$ commutes with the representatives of $\Aa$ and anticommutes with the Dirac operator.

\vspace{.1cm}

\noindent The predominant view of this paper is that Cantor sets should be treated as the boundary of a tree.  Using an idea of Michon \cite{MIC}, an ultrametric Cantor set $(C,d)$ defines a {\em weighted}, rooted tree graph $\Tt= \Tt(C,d)$, with a boundary that is isometrically equivalent to $C$.  Moreover, there is a one-to-one correspondence between a certain class of weighted, rooted trees and ultrametrics on a Cantor set. The tree graph then allows to define a Hilbert space and a Dirac operator which in turn defines a spectral triple for which $\Aa$ is the space of Lipschitz continuous functions on $(C,d)$. Conversely, the ultrametric can be recovered from the spectral triple.

\noindent The $\zeta$-function $\zeta(s) =\TR{(|D|^{-s})}$ of the Dirac operator is shown to be holomorphic in the domain $\Re(s) >s_0$ where $s_0$ is the upper box dimension of $(C,d)$. Following Connes \cite{CON,CO88}, the map $\mu(f) = \lim_{s\downarrow s_0} \TR{(f|D|^{-s})}/\TR{(|D|^{-s})}$ defines a probability measure on $C$. When $C$ is the attractor of a self-similar iterated function system system, this measure can be identified with the Hausdorff measure corresponding to the box dimension $s_0$. On the new Hilbert space $L^2(C,d\mu)$ there is an operator $\Delta_s$, $s\in\RM$, defined via the Dirichlet form $\langle-\Delta_sf,g\rangle = 1/2 \int_{\Upsilon(C)} \TR{(|D|^{-s}[D,\pi_\tau(f)]^\ast [D,\pi_\tau(g)])}d\nu(\tau)$.  Here $\Upsilon(C)$ is the space of {\em choices} on $C$ and serves as the analogue for $C$ of the unit sphere bundle of a manifold. It will be shown that for all $s\in\RM$, $\Delta_s$ is the generator of a Markovian semigroup on $L^2(C,d\mu)$ which can be seen as the analogue of a Brownian motion on $(C,d)$. 

\noindent The present work grew out of the authors' desire to create a spectral triple for the transversal, $\Xi$, of an aperiodic Delone set of finite type \cite{BBG, Bel, BHZ}.  In this case, the transversal is an ultrametric Cantor set, but it is not clear that $\Xi$ should be embeddable in $\RM^d$ for any $d\in \NM^+$.  Therefore, it was necessary to create a spectral triple that was reliant on the intrinsic data of the Cantor set.  The original idea came from a proposal by Alain Connes \cite{CON}(Chap. 4.3.$\epsilon$) for the triadic Cantor set. He also used this formalism to describe some properties of the Julia set \cite{CON}. In \cite{LAP}, Michel Lapidus proposed a program for applying the techniques of noncommutative geometry to fractals.  More specifically, he was interested in creating spectral triples on fractals that would recapture both the geometric properties (i.e. fractal dimension, Hausdorff measure, etc. \cite{FAL}) and the analytic properties(i.e. Laplacian \cite{KIG}) of the fractal. Since Lapidus outlined his program, there have been spectral triples proposed for many different types of fractals.  In \cite{CIL}, Lapidus and his coauthors were able to perform much of his program for the Sierpinski Gasket.  For the Cantor set, most of the work in this direction has based itself on the spectral triple proposed by Connes  for the triadic Cantor set.  Connes' spectral triple as well as many others based on his (e.g. \cite{GuI1},\cite{GuI2}) require the Cantor set to be embedded in $\RM^d$ for some $d\in\NM^+$.  In another direction, Christensen and Ivan  \cite{ChI} have proposed spectral triples for abstract compact metric spaces by gluing together spectral triples associated with pairs of points.  They are able to recover an equivalent metric with one of their spectral triples, but are unable to recover precise geometric information about their original space.  The spectral triple given in the present paper extends much of the previous work by recovering the appropriate fractal geometric information intrinsically without the necessity of embedding in $\RM^d$.  More importantly, this spectral triple allows to construct an analogue of the Laplacian on the Cantor set and therefore an analogue of Brownian motion.  

\noindent Diffusion on Cantor sets is not entirely new and has been studied in various contexts mostly as a non-Archimedean field \cite{AlK,Eva,Koc}.    Del Muto and Fig\`a-Talamanca have generalized this in \cite{dMF,Fig} for locally compact ultrametric spaces where the group of isometries is transitive and therefore allows to treat the Cantor set as an abelian group.  In both cases, the construction of the diffusion relies heavily on the algebraic structure that is given to the space.  A particular interest in these constructions has been taken by physicists for its potential applications in creating a $p$-adic based model of space-time based on $p$-adic differential operators \cite{VVZ}.  One of the most used and simplest of these operators is the Vladimirov operator \cite{VVZ} which has been used as the analogue of the Laplacian.  In the present paper, it is shown that the algebraic structure is unnecessary for the construction of an appropriate Laplacian.  In particular, the transversal of the Fibonacci tiling has only one nontrivial isometry and therefore has no obvious algebraic structure.  In the treatment of the example of the triadic Cantor set in the present paper, it is shown that the Vladimirov operator is related to the phase of the Dirac operator and consequently forgets the information provided by the metric.  Moreover, an eigenbasis for $\Delta_s$ is constructed and shown to be the basis of Haar wavelets thus recovering the result of \cite{Koz} for the Vladimirov operator.  In another area, Favre and Jonsson \cite{FaJ} have used an ultrametric tree to analyze singularities of algebraic varieties.  They define a Laplacian but it is unclear the relation to the present construction. The present paper then extends the work to date by providing precise asymptotic estimates as $t\to 0$ for the Brownian motion on the Cantor set.  More examples, including the case of the transversal of an aperiodic, repetitive Delone set of finite type, will be studied in future papers. 

\vspace{.5cm}
\section{Statement of Main Results}
\label{cantor07.sect-main}

\noindent This section presents a brief summary of the main results of this paper.  As stated in the introduction, the main viewpoint of this paper is that Cantor sets should be treated as the boundary of a tree.

\vspace{.1cm}

\noindent A Cantor set $C$ is a totally disconnected, compact, metrizable space without isolated points. It is well-known that such a set is homeomorphic to $C_0 =\{0,1\}^\NM$. Therefore, up to homeomorphism, the Cantor set is unique. However, adding a metric changes this prospect entirely.  The structure of the tree allows to capture the additional information provided by the metric.

\begin{defini}
\label{cantor07.def-sum}
Let $C$ be a Cantor set. A metric on $C$ will be called {\em regular} if it defines a topology on $C$ for which $C$ is a Cantor set.  A metric on $C$ is an {\em ultrametric} if $d(x,y) \leq \max\{d(x,z),d(z,y)\}$ for all $x,y,z \in C$.
\end{defini}

\noindent  Let now $(C,d)$ be a metric Cantor set such that $d$ is a regular ultrametric.  As mentioned previously, Michon was able to show that every ultrametric gives a weighted, rooted tree such that the boundary of the tree, $\pT$, is isometric to $C$.  The set of vertices and edges of its Michon tree $\Tt$ will be denoted by $\vs,\es$.  Any $v\in\Vv$ defines a clopen set $[v]\subset \pT$ which is the set of all infinite paths starting at the root that contain $v$.  It is then possible to build an even spectral triple $(\Aa,\Hh,D)$. The $\ast$-algebra $\Aa$ will be chosen as the space $\Aa= \lip(C) $ of Lipschitz continuous complex valued functions defined on $C$. From the noncommutative standpoint, it is important to note that this algebra is dense in the \Cs $\Cc(C)$ of continuous functions on $C$ and is invariant by the holomorphic functional calculus \cite{Bla86}.  In particular, this implies that the $K$-theory of $\Aa$ is the same as $\Cc(C)$. The Hilbert space $\Hh$ is given by $\Hh= \ell^2(\Vv)\otimes \CM^2$ where $\Vv$ is the set of vertices of $\Tt$. The grading operator is the multiplication by $\Gamma=\id \otimes \sigma_3$ where $\sigma_3 = \mbox{\rm diag}\{+1,-1\}$ is the third Pauli matrix. The Dirac operator $D$ is the operator defined by $D\psi\;(v) = (\diam [v])^{-1} \sigma_1 \psi(v)$ where $\sigma_1$ is the first Pauli matrix (equal to $+1$ off the diagonal and $0$ on the diagonal). To define the representation of the algebra $\Aa$ the notion of {\em choice} is needed.

\begin{defini}
\label{cantor07.def-choice}
Let $C$ be a Cantor set with a regular ultrametric $d$. A choice function is a map $\tau: \vs\mapsto C\times C$ such that, if $v\in \vs$ and if $\tau(v) = (x,y)$, then both $x,y$ are in $[v]$ and $d(x,y)= \diam [v]$. The set of choice functions on $C$ will be denoted by $\Upsilon(C)$. 
\end{defini}

\noindent In what follows $\tau(v)=(x,y)$ will be written $x=\tau_+(v),y=\tau_-(v)$. Then the $\ast$-representation $\pi_\tau$ of $\Aa$ is given by $\pi_\tau(f) \psi\;(v) = \mbox{\rm diag}\big\{f(\tau_+(v))\,,\,f(\tau_-(v))\big\}\; \psi(v)$.  It is important to note that the Dirac operator is independent of the choice function.  In the Noncommutative Riemannian structure, the space of choices $\Upsilon(C)$ will play the role of the unit sphere subbundle of the tangent bundle.  In particular the basic element of intuition is that a choice function is the analogue of a vector field of unit vectors on a manifold.  With this intuition in mind, then $[D,\pi_\tau(f)]$ represents the directional derivative of $f$ in the direction of $\tau$.  On $\RM^d$, a function $f\in\Cc^\infty(\RM^d)$ is such that the gradient $\|\nabla f\|_\infty < 1$ if and only if the directional derivative $\|\partial_{\vec{v}} f\|_\infty < 1$ for every $\vec{v}\in\RM$.  Therefore, by this reasoning it is natural to expect that the metric can be recovered by using the Connes distance on functions $f\in\Cc(C)$ such that $\|[D,\pi_\tau(f)]\| < 1$ for every $\tau\in\Upsilon(C)$.  This is shown in the following:

\begin{theo}
\label{cantor07.th-connesdis}
Let $C$ be a Cantor set with a regular ultrametric $d$. Then $d$ coincides with the Connes distance $\rho$ defined by 

$$\rho(x,y):=
   \sup\{
     |f(x)-f(y)| \,:\, f\in\lip(C)\,,\,
    \sup_{\tau\in\Upsilon(C)}
      \|[D,\pi_{\tau}(f)]\| \,\leq \, 1
      \}
$$
\end{theo}

\noindent This result indicates that the spectral triple defined above is sufficient to recover the metric $d$ on $C$ whenever $d$ is a regular ultrametric and when all possible choice functions are taken into account. In fact, $\rho$ is the typical Connes distance with respect to the spectral triple obtained by summing over all possible choice functions.

\vspace{.1cm}

\noindent 
Following the idea of Connes \cite{CO88,CON}, let $\zeta(s):=\TR(|D|^{-s})$ be the $\zeta$-function associated with the Dirac operator. It is known that there is $s_0 >0$ (possibly infinite), such that $\zeta$ is holomorphic with respect to $s$ in a half-plane of the form $\Re (s) >s_0$ and that $\zeta$ is singular at $s_0$. Then $s_0$ is called its {\em abscissa of convergence}. 

\vspace{.1cm}

\noindent Let $\Tt$ be the tree corresponding to $(C,d)$.  Let $\{\lambda_k\}_{k=1}^\infty$ be the set of all distinct $\tdiam([v])$ for $v\in\Vv$ (these are also the distinct eigenvalues of $|D|^{-1}$).  Let them be ordered such that $\lambda_1 > \lambda_2 >\cdots$.  Let $M_n$ be such that every vertex with diameter at least $\lambda_n$ has at most $M_n$ children.  Then the next result is the following:

\begin{theo}
\label{cantor07.th-box}
If $(\log M_n)/(-\log \lambda_n)\to 0$ as $n\to \infty$, then $s_0 = \ubox(C)$.
\end{theo}

\noindent It is important to note that a special case of Theorem 2 is when there is a uniform bound on the number of children - this happens for the attractor of a self-similar iterated function system and the transversal of the Fibonacci tiling.  In any case, the hypothesis says intuitively that the number of children can grow but it must be compensated for by a  decrease in the size of the children.

\vspace{.1cm}

\noindent At last, $(C,d)$ will be called {\em $\zeta$-regular} whenever the abscissa of convergence of the $\zeta$-function is finite, if $\lim_{s\downarrow s_0}(s-s_0)\TR\left( |D|^{-s} \right) > 0$ and if, for any $f\in \Aa$ the following limit exists

\begin{equation}
\label{cantor07.eq-mu}
\mu(f) \;=\;
   \lim_{s\downarrow s_0}
    \frac{\TR\left( |D|^{-s} \pi_\tau(f)\right)}
     {\TR\left( |D|^{-s} \right)}
\end{equation}
\begin{theo}
\label{cantor07.th-meas}
Let $C$ be a $\zeta$-regular Cantor set with a regular ultrametric $d$. Then  the limit (\ref{cantor07.eq-mu}) is independent of the choice function $\tau$ and defines a probability measure on $C$. 
\end{theo}

\noindent  Given the measure $\mu$ it is possible to construct various operators on $L^2(C, \mu)$.  In order to do so, it is necessary to use $\mu$ to define a measure $\nu$ on the space of choices $\Upsilon(C)$.  It is then interesting to try to find the analogue of the Laplace-Beltrami operator on a Riemannian manifold.  Locally, the Laplace-Beltrami operator on a Riemannian manifold is the average of the square of the directional derivatives, $|\nabla_{\vec v}f|^2$, over the unit sphere of the tangent space.  However, the Euclidean structure of the tangent space allows this average to be reduced to a sum of the directional derivates in the direction of an orthonormal basis.  In the case of the Cantor set, the local basis is infinite and is given by choice functions.  Therefore since there is no local Euclidean structure on $\Upsilon(C)$, it is natural to define the analogue of the Laplace-Beltrami operator by the following:  

\begin{theo}
\label{cantor07.th-lapla}
Let $C$ be a $\zeta$-regular Cantor set with a regular ultrametric $d$.  Then the measure $\mu$ coming from the $\zeta$-function defines a measure $\nu$ on the space of choices $\Upsilon(C)$.  Moreover, for all $s\in\RM$ there is a closable Dirichlet form on the Hilbert space $L^2(C, \mu)$ defined by  

$$Q_s(f,g):=\frac{1}{2}\int_{\Upsilon(C)}\TR(|D|^{-s}[D,\pi_\tau(f)]^*[D,\pi_\tau(g)])d\nu(\tau)$$

\noindent with $\Dom(Q_s)$ a dense subspace of the real Hilbert space $L^2(C,\mu)$.
 
\end{theo}

\noindent When $\Supp(\mu)=C$ then the classical theory of Dirichlet forms \cite{Fuk} associates to $Q_s$ a non-positive definite self-adjoint operator $\Delta_s$ on $L^2(C,\mu)$ which generates a Markovian semigroup.  It will be shown in the case of the triadic Cantor set $C_3$ that for $s=s_0$, $\Delta_{s_0}$ plays the role of a Laplacian on $C_3$ in the sense that the Weyl asymptotic formula gives $\Nn(\lambda)\sim c_0\lambda^{s_0/2}$.   The Markovian semigroup associated to $\Delta_{s_0}$ defines a stochastic process $\big(X_t\big)_{t\geq 0}$ with values in $C_3$.  However, somewhat unexpectedly this Brownian motion on $C_3$ satisfies $\EM\big(d(X_{t_0},X_{t_0+t})^2\big) \sim c_1 t\ln(1/t)$ as $t\downarrow 0$.  This surprising subdominant contribution by $\ln(1/t)$ needs further investigation.

\vspace{.5cm}
\section{Rooted Trees}
\label{cantor07.sect-trees}

 \subsection{Basic Definitions}
 \label{cantor07.ssect-bdef}

\noindent This section is a reminder about rooted trees (see \cite{BOL}). A graph is a triple $G=\{\vs,\es,\psi\}$ where $\vs$ is a non-empty countable set with elements called {\em vertices} and $\es$ is a countable set with elements called {\em edges}. Let $\vs^{(2)}$ denote the set of unordered pairs of vertices, namely $\vs^{(2)} = \vs\times \vs/\sim$ where $\sim$ is the equivalence relation defined by $(v,w)\sim (w,v)$. Then $\psi:\es\to \vs^{(2)}$ is called the {\em incidence function} which assigns to each edge an unordered pair of not necessarily distinct vertices. If $\psi(e) = (v,w)$ then $e$ is said to {\em link} $v$ and $w$, while $\tilde{\psi}(e)$ will denote the set $\{v,w\}$. The {\em degree} $|v|$ of a vertex $v\in \vs$ is the number of edges $e\in \es$ such that $v\in\tilde{\psi}(e)$. Two vertices $v,w\in \vs$ will be called {\em incident} if there exists an edge $e\in \es$ such that $\psi(e)=(v,w)$. A graph $G$ is {\em simple} whenever (i) there are no edges $e$ such that $\psi(e) = (v,v)$ for some $v\in \vs$ and (ii) if $e,e'$ are two edges such that $\psi (e)= \psi(e')$, then $e=e'$. In what follows only simple graphs will be considered.

\vspace{.1cm}

\noindent A {\em walk} on a graph $G$ is a double sequence $\{(v_0, v_1, \cdots, v_{n-1}, v_n)\,;\, (e_1,e_2,\cdots ,e_n)\}$ (where $n$ is finite or infinite) of incident vertices and edges linking them  such that $\psi(e_i) =(v_{i-1},v_i)$ for all $i>0$. For a simple graph it is sufficient to specify the sequence of vertices. A {\em step} of the walk is a triple of the form $(v_{i-1}, v_i, e_i)$. The length of the walk is the number $n$ of steps making this walk. If the walk is finite, the first and the last vertices of the sequence are said to be {\em linked} by the walk. If $v$ is one of the vertices of the walk, then the latter is said to {\em pass} through $v$. A {\em path} is a walk with pairwise distinct vertices. The graph $G$ is {\em connected} if given any two vertices there is a finite path linking them.

\vspace{.1cm}

\noindent A {\em cycle} is a finite walk with at least three steps, such that the first and the last vertices coincide and all other vertices are pairwise distinct. A {\em tree} is a connected graph with no cycle. A {\em rooted tree} is a pair $(\Tt, 0)$ where $\Tt$ is a tree and $0$ is a vertex of $\Tt$ called the {\em root}. By abuse of notation, $\Tt$ will denote a rooted tree, and the root will be implicit. Since $\Tt$ is a tree, given any pair of distinct vertices there is one and only one path linking them. In particular, there is a unique path linking the root to a given vertex. So that there is a one-to-one correspondence between the set of vertices and the set of finite paths starting at the root. 

\vspace{.1cm}

\noindent On a rooted tree there is a partial order defined by $v\succeq w$ if the path from the root to $w$ necessarily passes through $v$. Then $w$ is called a {\em descendant} of $v$ and this will also be written as $w\preceq v$, while $v$ will be called an {\em ancestor} of $w$. If, in addition, $v,w$ are incident, then $v$ is called the {\em father} of $w$ and $w$ is called a {\em child} of $v$. The {\em height}, $\theight(v)$, of a vertex $v$ is the length of the unique path linking the root to $v$. Hence the root has height $0$, its children have height $1$ and so on.

\vspace{.3cm}  

 \subsection{The Boundary of a Rooted Tree}
 \label{cantor07.ssect-bound}

\noindent In this section $\Tt$ will denote an {\em infinite} rooted tree with root $0$. The set $\vs$ of its vertices is endowed with the discrete topology. Since it is infinite it is certainly not compact. A compactification of the tree can be defined by considering the {\em boundary} $\pT$ of this tree defined as follows:

\begin{defini}
\label{cantor07.def-btree}
If $\Tt$ is a rooted tree, its boundary $\pT$ is the set of infinite paths starting at the root. 
\end{defini}

\noindent A vertex is {\em dangling} if it has no child. Hence the boundary ignores dangling vertices. In what follows, only trees with no dangling vertices will be considered. This implies among other things that every finite path can be extended to an infinite path.

\begin{exam}
\label{cantor07.exam-bin}
{\em Let $T_2$ be the infinite binary rooted tree.  That is $T_2$ is the tree with a root and such that every vertex has exactly two children.  Since every vertex has two children the edge linking it to one child will labeled by $0$ and the other by $1$. Hence any finite path starting at the root, and therefore any vertex, is labeled by a finite sequence of $0$'s and $1$'s. The root is given by the empty sequence. Thus $\Tt_2$ can be seen as the set of finite sequences of of $0$'s and $1$'s. Consequently, $\pT_2=\{0,1\}^\NM$. The map $\Theta : (\epsilon_n)_{n\in \NM} \mapsto \sum_{i=0}^\infty 2\epsilon_i\,3^{-(i+1)}$ defines a one-to-one map from $\pT_2$ onto the classical triadic Cantor set.
}
\hfill $\Box$
\end{exam}

\begin{defini}
\label{cantor07.def-clop}
Let $\Tt$ be a rooted tree. If $v$ is a vertex, $[v]\subset \pT$ denotes the set of infinite paths starting at the root and passing through $v$. 
\end{defini}

\begin{proposi}
\label{cantor07.prop-btree}
Let $\Tt$ be a rooted tree with no dangling vertex. Then, the set $\{[v]\,;\, v\in \vs\}$ is a basis of open sets for a topology on the boundary of $\Tt$ for which $\pT$ is completely disconnected. For this topology $\pT$ is compact if and only if  each vertex has at most a finite number of children. It has no isolated points if and only if each vertex has one descendant with at least two children.
\end{proposi}

\noindent {\bf Proof: } (i) Clearly the family covers $\pT$ since $[0]=\pT$. Moreover, if $v,w\in\vs$ the intersection $[v]\cap [w]$ is either empty or, if not, then one of the two vertices is an ancestor of the other. In particular if, say $v\succeq w$, then $[v]\cap [w]= [w]$, showing that indeed this family is a basis for a topology on $\pT$. 

\noindent (ii) Let $v\in\vs$. Then let $\vs(v)$ be the set of vertices with same height as $v$. Clearly if $w\neq v$ and $w\in \vs(v)$ then $w$ is not comparable to $v$, hence $[v]\cap [w]=\emptyset$. Moreover, if $x\in \pT$ is an infinite path starting at the root, one of its vertices, say $w$, is such that $w\in\vs(v)$ and $x\in [w]$. Consequently, the family $\{[w]\,;\, w\in \vs(v)\}$ is a partition made of open sets. In particular, the complement of $[v]$ is the union of open sets and is open as well. Hence, for any vertex $v$, the set $[v]$ is a closed and open set (or a {\em clopen} set), so that $\pT$ is completely disconnected.

\noindent (iii) If there is a vertex $v$ having an infinite number of children, the family of $\{[w]\}$ such that $w$ is a child of $v$ defines an open covering of $[v]$ from which no finite covering can be extracted since this is a partition. Thus $[v]$, which is closed, cannot be compact and thus $\pT$ cannot be compact either. 

\noindent (iv) Conversely, let $\Tt$ be such that each of its vertices has only finitely many children and let $\Oo$ be an open cover of $\pT$.  There exists an $N$ such that, for each $v\in \Vv$ of height $N$, there is an $O_v\in\Oo$ with $[v]\subset O_v$. Suppose not.  Then there exists a sequence of vertices $v_0v_1\cdots$ such that each $[v_k]$ is not covered by any $O\in\Oo$. Moreover, this sequence actually gives an infinite path $\sigma=v'_0v'_1\cdots$ such that each $v'_k$ is not covered by any single $O\in\Oo$.  This path is constructed as follows.  One of the children of the root, called $v'_1$, must contain an infinite number of $v_k$.  In the same way, one of the children of $v'_1$, called $v'_2$ must contain an infinite number of $v_k$.  Proceeding recursively, an infinite sequence $v'_0v'_1\cdots$ is obtained such that (i) for each $n\geq 0$, $v'_k$ is a child of $v'_{k-1}$ and (ii) $[v'_k]$ is not covered by any single $O\in\Oo$.  Then $v'_0v'_1\cdots\in\pT$ and is not covered by $\Oo$ which contradicts the fact that $\Oo$ is an open cover.  Consequently, since each vertex has only a finite number of children, then there are only a finite number of vertices of height $N$.  Therefore, $\Oo$ has a finite subcover and $\pT$ is compact.

\noindent (v) Let $v$ be a vertex of $\Tt$ such that none of its descendants has more than one child. Then $[v]$ is reduced to one single path $x$ which is itself an open set. Hence $x$ is isolated. Conversely, if $x\in\pT$ is isolated, then $\{x\}$ is open, meaning that it contains at least one nonempty element of the basis. Hence there is $v\in \vs$ such that $[v] \subset\{x\}$. But this can happen only if each descendant of $v$ has only one child, since otherwise, $[v]$ would contain at least two distinct infinite paths. 
\hfill $\Box$

\begin{defini}
\label{cantor07.def-cantree}
A tree will be called Cantorian if it has a root, no dangling vertex and if each vertex has a finite number of children as well as a descendant with more than one child.
\end{defini}

\begin{rem}
\label{cantor07.rem-cantree}
{\em By Prop.~\ref{cantor07.prop-btree} this definition is equivalent to $\pT$ is a Cantor set.
}
\hfill $\Box$
\end{rem}

\noindent Various surgical operations on a tree lead to similar boundaries. The first operation is {\em edge reduction}. Namely if there is a path $\gamma$ linking $v$ to one of its descendant $w$ such that each vertex of this path distinct from $v,w$ has only one child, then the graph can be reduced by suppressing these vertices and replacing the path by one edge. Hence if $x\in\pT$ is any path passing through $v$ and $w$, it also automatically passes through all of the vertices of $\gamma$. Then it can also be reduced and the reduction operation gives a one-to-one mapping between the boundary of the initial tree and the boundary of the reduced one. In addition $[u]=[w]$ whenever $v\succeq u \succeq w$ and $u\neq w$, so that this mapping is actually an homeomorphism. 

\noindent The opposite of edge reduction will be called {\em edge extension}. Namely any edge can be replaced by a finite path with same end points so that each internal vertex of the path has only one child. 

\noindent There is also the notion of {\em vertex extension}. Namely if $v$ is a vertex with at least three children then one child will be called $v_0$ and the others $v_1,\cdots, v_r$. Then a new vertex $u$ is created as a child of $v$ having $v_1,\cdots, v_r$ as children. As before, this vertex extension does define also an homeomorphism between the corresponding boundaries. In particular this implies the following proposition which is one of many ways of showing that every Cantor set is homeomorphic to $\{0,1\}^{\NM}$.

\begin{proposi}
\label{cantor07.prop-bin}
Let $\Tt$ be a Cantorian tree. Then there is a map made of the product of a possibly infinite family of edge reductions, edge extensions and vertex extensions, mapping $\Tt$ onto the binary tree $\Tt_2$ and defining a homeomorphism of their boundaries. 
\end{proposi}

\begin{defini}
\label{cantor07.def-lcp}
Let $\Tt$ be a Cantorian tree. If $A\subset \pT$ then a vertex $v$ is a common ancestor of $A$ if $A\subset [v]$. If $A$ has more than one point, its least common prefix (or l.c.p.) is the smallest of its ancestors. If $A=\{x,y\}$ the least common prefix will be denoted by $x\wedge y$. 
\end{defini}

\begin{proposi}
\label{cantor07.prop-lcp}
Let $\Tt$ be a Cantorian tree. The l.c.p. of a subset $A\subset \pT$ with more than one point always exists and is unique.
\end{proposi}

\noindent {\bf Proof: } Since $[0]=\pT$ it follows that $A$ always admits the root as an ancestor. Now if $v$ and $w$ are both common ancestors of $A$, then since $A \subset [v]\cap [w]$ is non empty it follows that one of the two vertices, say $v$ is an ancestor of the other, so that $A \subset [w] \subset [v]$. Hence the set of common ancestors of $A$ is totally ordered. Since it is at most countable this set defines a path $x=(0=v_0,v_1,\cdots ,v_n)$. Since $A$ contains at least two distinct points, this path is automatically finite because otherwise the intersection $\bigcap_{i\geq 0} [v_i]$ would be reduced to $\{x\}$ and would contain $A$, a contradiction. Thus $v_n$ is the least common ancestor and is unique. 
\hfill $\Box$

\vspace{.5cm}
\section{Michon's Correspondence}
\label{cantor07.sect-mc}

\noindent For the sake of the reader, this section recalls Michon's correspondence between regular ultrametrics on a Cantor set $C$, profinite structures on $C$, and weighted, rooted trees.

\vspace{.3cm}  

 \subsection{Ultrametrics and Profinite Structures}
 \label{cantor07.ssect-ups}

\noindent This section shows the correspondence between ultrametrics and profinite structures on $C$ \cite{MIC}
.  Let $C$ be a Cantor set with regular metric $d$.  Following \cite{HAU}, given $\epsilon > 0$ and $x,y\in C$ let an {\em $\epsilon$-chain} be a sequence $x_0=x, x_1,\dots x_{n-1}, x_n=y$ of points in $C$ such that $d(x_i,x_{i+1}) < \epsilon$.  This gives rise to an equivalence relation $\stackrel{\epsilon}{\sim}$ by defining $x\stackrel{\epsilon}{\sim} y$ if there is an $\epsilon$-chain between them.  In such a case, $[x]_\epsilon$ will denote the equivalence class of $x\in C$.  It is then possible to define the {\em separation} of $x$ and $y$ by $\delta(x,y):=\inf\{\epsilon: x\stackrel{\epsilon}{\sim} y\}$.

\begin{proposi}
\label{cantor07.prop-sep}
Let $C$ be a Cantor set with regular metric $d$. Then the separation $\delta$ is the maximum ultrametric on $C$ dominated by $d$.  Moreover, $\delta$ is regular.
\end{proposi}

\noindent {\bf Proof: } 
By \cite{HAU} (Ch 29.3), $\delta$ is an ultrametric on the connected components.  Since $C$ is totally disconnected then $\delta$ is an ultrametric on $C$.  If $d(x,y)=\epsilon$ then $x \stackrel{\epsilon}{\sim} y$.  Therefore, $\delta(x,y) \leq d(x,y)$.  Now let $d'$ be another ultrametric on $C$ such that $d'(x,y)\leq d(x,y)$ for $x,y\in C$.  Then for any $\epsilon$-chain $x_0=x,\dots,x_n=y$, $$d'(x,y)\leq\max\{d'(x_i,x_{i+1}): 0\leq i\leq n-1\}\leq \max\{d(x_i,x_{i+1}): 0\leq i\leq n-1\}<\epsilon.$$  Thus, $d'\leq \delta$.  For a proof that $\delta$ is regular see \cite{HAU}.
\hfill $\Box$

\vspace{.3 cm}

\noindent  It follows at once from the proposition that if $d$ is an ultrametric then $d=\delta$.  From now on, it will be assumed that $C$ is a Cantor set with regular ultrametric $d$.

\begin{defini}
\label{cantor07.def-pfs}
A {\em profinite structure} on a Cantor set $C$ is given by an increasing family $\{R_\epsilon:\epsilon\in\RM^+\}$ of equivalence relations on $C$ that satisfy the following properties:

\noindent(i) Each relation $R_\epsilon$ is open in $C\times C$ and for a certain $\epsilon$, $R_\epsilon=C\times C$;

\noindent(ii) The family is continuous on the left: $\bigcup_{\epsilon' < \epsilon} R_{\epsilon'}=R_\epsilon$;

\noindent(iii) $\bigcap_{\epsilon\in\RM^+} R_\epsilon = \Delta$ (the diagonal of $C\times C$).
\end{defini}

\begin{proposi}
\label{cantor07.prop-umpfs}
On a Cantor set $C$, there is a one-to-one correspondence between profinite structures and regular ultrametrics.
\end{proposi}

\noindent A proof of this result is given in the appendix.

\vspace{.3 cm}

 \subsection{Weighted, Rooted Trees}
 \label{cantor07.ssect-wrt}

Using the results of the last section, it is now possible to show the connection between Cantorian trees and ultrametrics on a Cantor set.

\begin{defini}
\label{cantor07.def-wt}
Let $\Tt$ be an infinite rooted tree with no dangling vertex.  A {\em weight} on $\Tt$ is a function $\epsilon:\Vv\to\RM^+$ that satisfies the following:

\noindent(i) If $v\succ v'$ then $\epsilon(v) > \epsilon(v')$.

\noindent(ii) For an infinite path $v_0v_1\cdots\in\pT$, $\lim_{n\to\infty} \epsilon(v)=0$.

\noindent A rooted tree along with its weight function will be called a {\em weighted}, rooted tree.
\end{defini}

\noindent As mentioned previously, there are various surgical operations on trees that lead to the same boundary.  Given a tree $\Tt$, any vertex with only one child can be reduced by the process of edge reduction.  The weight function is then the restriction of the original weight function.  A tree for which every vertex has at least two children will be called {\em reduced}.  

\begin{proposi}
\label{cantor07.prop-wrtpfs}
On a Cantor set $C$, there is a one-to-one correspondence between regular ultrametrics and reduced, weighted, rooted Cantorian trees.  Moreover given a regular ultrametric $d$, the boundary $\pT$ of the corresponding weighted, rooted Cantorian tree is isometric to $(C,d)$.  The weight function $\epsilon$ for $\Tt$ is such that $\epsilon(v)=\tdiam_d([v])$.
\end{proposi}

\noindent A proof of this result is given in the appendix.

\vspace{.3cm}

\subsection{Embedding of Ultrametric Cantor Sets}
\label{cantor07.ssect-embedding}

\noindent  A simple application of Michon's correspondence is given by the following.

\begin{theo}
\label{cantor07.th-embedding}
Let $C$ be a Cantor set with regular ultrametric $d$.  Let $\Tt$ with weight $\epsilon$ be the corresponding reduced, weighted, rooted Cantorian tree.  If $\Vv_*$ denotes all the vertices of $\Tt$ except for the root, then there exists an isometric embedding of $C$ into the real Hilbert space $\ell^2_\RM(\Vv^*)$.
\end{theo}

\noindent {\bf Proof: } Let $x\in C$ and let $v_0v_1\cdots$ be the infinite path corresponding to $x$.  Let
$$\Phi(x):=\sum_{n=0}^\infty\sqrt{\frac{\epsilon(v_n)^2-\epsilon(v_{n+1})^2}{2}}|v_{n+1}\rangle$$
where $\{|v\rangle, v\in\Vv_*\}$ denotes the canonical basis of $\ell^2_\RM(\Vv^*)$.  If $v\neq v'$, then $\langle v,v'\rangle = 0$.  Therefore,
$$||\Phi(x)||^2=\sum_{n=0}^\infty\frac{\epsilon(v_n)^2-\epsilon(v_{n+1})^2}{2}=\frac{\epsilon(v_0)^2}{2}$$
and $\Phi(x)\in \ell^2_\RM(\Vv^*)$.  Thus, $\Phi$ is well-defined.  Let $y\in C$ with $y\neq x$.  If $w_0w_1\cdots$ is the infinite path corresponding to $y$ then there exists an $n_0 > 0$ such that $w_n\neq v_n$ for $n > n_0$ and $w_n=v_n$ for $n \leq n_0$.  Then $x\wedge y= v_{n_0}$ and $d(x,y)=\epsilon(v_{n_0})^2$.  Moreover,

$$\Phi(x)-\Phi(y)=
   \sum_{n=n_0}^\infty
    \sqrt{\frac{\epsilon(v_n)^2-\epsilon(v_{n+1})^2}{2}} |v_{n+1} \rangle
    -\sum_{n=n_0}^\infty\sqrt{\frac{\epsilon(w_n)^2-\epsilon(w_{n+1})^2}{2}} |w_{n+1} \rangle
$$

\noindent and consequently

$$\|\Phi(x)-\Phi(y)\|^2=\sum_{n=n_0}^\infty\frac{\epsilon(v_n)^2-\epsilon(v_{n+1})^2}{2}+\sum_{n=n_0}^\infty\frac{\epsilon(w_n)^2-\epsilon(w_{n+1})^2}{2}=\epsilon(v_{n_0})^2.
$$

\noindent Since $\epsilon(v_{n_0})=d(x,y)$ then $\Phi$ is indeed an isometry.
\hfill $\Box$

\vspace{.5cm}
\section{A Spectral Triple}
\label{cantor07.sect-spt}

\noindent Given Michon's correspondence, it is now possible to construct a spectral triple on a Cantor set $C$ with regular ultrametric $d$.

\vspace{.3cm} 

 \subsection{Construction of the Spectral Triple}
 \label{cantor07.ssect-spectint}

\begin{defini}
\label{cantor07.def-st}
An {\em odd spectral triple} for an involutive algebra $\Aa$ is a triple $(\Aa,\Hh,D)$ where $\Hh$ is a Hilbert space on which $\Aa$ has a representation $\pi$ by bounded operators.  $D$ is a self-adjoint operator on $\Hh$ such that $[D,\pi(a)]$ is a bounded operator on $\Hh$ for all $a\in\Aa$ and such that $(D^2+1)^{-1}$  is compact.

An {\em even spectral triple} is an odd spectral triple along with a grading operator $\Gamma:\Hh\to\Hh$.  $\Gamma$ is required to satisfy $\Gamma^*=\Gamma$, $\Gamma^2=1$, $\Gamma D=-D\Gamma$, and $\Gamma\pi(a)=\pi(a)\Gamma$ for all $a\in\Aa$
\end{defini}

\noindent The algebra will be $\lip(C)$.  Let $\Tt$ be the reduced, weighted, rooted Cantorian tree corresponding to the regular ultrametric $d$.  Since $\Tt$ is Cantorian, the set of vertices $\Vv$ is countable.  Let $\Hh:=\ell^2(\Vv)\otimes\CM^2$.  $D$ is the operator on $\Hh$ given by $D\psi\;(v):=(\tdiam(v))^{-1}\sigma_1\psi(v)$ where $\sigma_1$ is the first Pauli matrix.  The grading operator is the multiplication by $\Gamma:=\id \otimes \sigma_3$ where $\sigma_3 = \mbox{\rm diag}\{+1,-1\}$ is the third Pauli matrix.  As mentioned earlier, to define a representation on $\Aa$ a notion of choice is required.  Let $\tau\in\Upsilon(C)$ be a choice function.  Then the $\ast$-representation $\pi_\tau$ of $\Aa$ is given by $\pi_\tau(f) \psi\;(v) = \mbox{\rm diag}\big\{f(\tau_+(v))\,,\,f(\tau_-(v))\big\}\; \psi(v)$.

\begin{proposi}
\label{cantor07.prop-fr}
$\pi_\tau$ is a faithful $\ast$-representation of $\Cc(C)$ for all $\tau\in\Upsilon(C)$.
\end{proposi}

\noindent {\bf Proof: } That $\pi_\tau$ is a $\ast$-representation is obvious.  It is bounded since $f$ is continuous and $C$ is compact.  Let $f,g\in\Cc(C)$ be such that $\pi_\tau(f)=\pi_\tau(g)$.  Then $f(\tau_+(v))=g(\tau_+(v))$ for all $v\in\Vv$.  For $x\in C$, there exists $v_0,v_1,\dots\in\Vv$ such that $x\in[v_j]$ and $\tdiam([v_j])\to 0$.  Then $f(x)=\lim_{j\to\infty} f(\tau_+(v_j))=\lim_{j\to\infty} g(\tau_+(v_j))=g(x)$.  Thus $\pi_{\tau}$ is faithful.
\hfill $\Box$

\vspace{.3 cm}

\noindent Based on this proposition, $\pi_\tau$ is also a faithful representation on $\lip(C)$.

\begin{proposi}
\label{cantor07.prop-spect}
$(\lip(C),\Hh,D,\Gamma)$ is an even spectral triple for all $\tau\in\Upsilon(C)$.
\end{proposi}

\noindent {\bf Proof: } To show that $D$ is self-adjoint, let $\psi,\psi'\in\Hh$.  Then,
$$\langle D\psi,\psi'\rangle_{\Hh}=\sum_{v\in \Vv} {(\tdiam([v]))^{-1}}\langle\sigma_1\psi(v),\psi'(v)\rangle_{\CM^2}=\langle \psi,D\psi'\rangle_{\Hh}$$
since $\sigma_1^*=\sigma_1$. Since $D$ is densely defined then $D$ is symmetric.  By \cite{CONW}(Prop. X.2.4), if Range$(D)=\Hh$ then $D$ is self-adjoint.  Let $\psi\in\Hh$ and let $\psi'(v)=\tdiam(v)\sigma_1\psi(v)$.  Then $D\psi'(v)=\psi(v)$ since $\sigma_1^2=\id$.  Now, since there exists $K$ such that $\tdiam(C)\leq K$, then
$$||\psi'||_\Hh^2=\sum_{v\in \Vv }(\tdiam([v]))^2||\psi(v)||_{\CM^2}\leq K^2||\psi||_\Hh^2$$
So, $\psi'\in\Hh$ and Range$(D)=\Hh$.  Thus $D$ is self-adjoint.

\noindent Let $v\in \Vv$.  Because $\tdiam([v])=d(\tau_+(v),\tau_-(v))$ since $\tau$ is a choice function, then
$$
([D,\pi_\tau(f)]\psi)(v)=\frac{f(\tau_+(v))-f(\tau_-(v))}{d(\tau_+(v),\tau_-(v))}\left(\begin{array}{ccc} 0 & -1 \\ 1 & 0 \end{array} \right)\psi(v).
$$ 
Since $f$ is Lipschitz, then $||[D,\pi_\tau(f)]\psi||_\Hh\leq k||\psi||_\Hh$ where $k$ is the Lipschitz constant of $f$ and $[D,\pi_\tau(f)]\in\Bb(\Hh)$.

\noindent To show that $(D^2+1)^{-1}$ is compact, let $\psi\in\Hh$ and $v\in \Vv$. Then, it is a straightforward calculation to show that
$$((D^2+1)^{-1}\psi)(v)=\frac{\tdiam([v])^2}{1+\tdiam([v])^2}\psi(v).$$
So for $\eta > 0$, let $(T^\eta\psi)(v)=((D^2+1)^{-1}\psi)(v)$ if $\tdiam([v]) \geq \sqrt\eta$ and $0$ otherwise.  Now since there are only finitely many $v\in \Vv$ with $\tdiam([v]) \geq \sqrt\eta$ then $T^\eta$ is finite rank.  Consequently,
$$
||T^\eta-(D^2+1)^{-1}||_{\Bb(\Hh)}=\sup_{v\in \Vv}\{\frac{\tdiam([v])^2}{1+\tdiam([v])^2}: \tdiam([v]) < \sqrt\eta\} < \eta
$$
Thus $||(D^2+1)^{-1}-T^\eta||_{\Bb(\Hh)} < \eta$ and $\lim_{\eta\downarrow 0} T^{\eta}=(D^2+1)^{-1}$.  Consequently $(D^2+1)^{-1}$ is compact. The proof that $\Gamma^*=\Gamma$, $\Gamma^2=\Gamma$, $\Gamma D=-D\Gamma$ and $\Gamma\pi_\tau(f)=\pi_\tau(f)\Gamma$ for all $f\in\lip(C)$ is straightforward.
\hfill $\Box$

\vspace{.3cm}

 \subsection{The Connes Distance: Proof of Theorem~\ref{cantor07.th-connesdis}}
 \label{cantor07.ssect-sub}

\noindent The spectral triple should be able to recover some of the structure of the original space $C$.  Theorem~\ref{cantor07.th-connesdis} shows that it can recover the metric when all possible choice functions are taken into account. 

\vspace{.3 cm}

\noindent {\bf Proof of Theorem~\ref{cantor07.th-connesdis}: }
Let $x,y\in C$ and let $d_x:C\to C$ be given by $d_x(y)=d(x,y)$.  Then $d_x$ is Lipschitz continuous.  Let $\tau\in\Upsilon(C)$ and recall that this implies that $d(\tau_+(v),\tau_-(v))=\tdiam(v)$.  Then 
$$
\|[D,\pi_{\tau}(d_x)]\|_{\Bb(\Hh)}=\sup_{v\in \Vv}\{\frac{|d(x,\tau_+(v))-d(x,\tau_-(v))|}{d(\tau_+(v),\tau_-(v))}\}
\leq\sup_{v\in \Vv}\{\frac{d(\tau_+(v),\tau_-(v))}{d(\tau_+(v),\tau_-(v))}\}=1
$$
where the inequality follows from the triangle inequality.  Consequently, $$\sup_{\tau\in\Upsilon(C)}\{||[D,\pi_{\tau}(d_x)]||_{\Bb(\Hh)}\}\leq 1$$ and $\rho(x,y)\geq |d_x(x)-d_x(y)|=d(x,y)$.

\noindent For $x,y\in C$, let $v\in \Vv$ be such that $v=x\wedge y$, so that $d(x,y)=\tdiam(v)$.  Let $\tau$ be such that $\tau_+(v)=x$ and $\tau_-(v)=y$.  Then for any $f\in\lip(C)$ such that $||[D,\pi_{\tau}(f)]||_{\Bb(\Hh)}\}\leq 1$
$$\frac{|f(\tau_+(v))-f(\tau_-(v))|}{\tdiam(v)}=\frac{|f(\tau_+(v))-f(\tau_-(v))|}{d(\tau_+(v),d(\tau_-(v))}\leq 1.$$
This gives that $|f(x)-f(y)|\leq d(x,y)$ and therefore that $\rho(x,y)\leq d(x,y)$.
\hfill $\Box$

\vspace{.5cm}
\section{$\zeta$-Functions}
\label{cantor07.sect-zeta}

\noindent In this section, the Dirac operator $D$ is used to create a $\zeta$-function as formulated by Connes \cite{CON}.  Since the Dirac operator is independent of choice, this $\zeta$-function will also be independent of choice.

\vspace{.3 cm}

 \subsection{The $\zeta$-function for $D$}
 \label{cantor07.ssect-zeta}

 Let $\Hh$ be the Hilbert space from the previously created spectral triple.  Then for $\psi\in\Hh$, $ (|D|\psi)(v)=\tdiam^{-1}(v)\psi(v)$.  Since $\tdiam(v) > 0$ for all $v\in \Vv$ then $|D|$ is invertible and $(|D|^{-1}\psi)(v)=\tdiam(v)\psi(v)$.  Let 
$$\zeta(s):=\frac{1}{2}\TR(|D|^{-s})=\sum_{v\in\Vv}\tdiam(v)^s.$$
Then $\zeta$ is a Dirichlet series.  By \cite{HaR} (Ch. 2), as a function of the complex variable $s$, $\zeta$ either converges everywhere, nowhere, or in a half-plane given by $\Re(s) > s_0$.  In the last case, $s_0$ is called the {\em abscissa of convergence}.  Since the eigenvalues of $|D|^{-1}$ are discrete, let $\zeta(s)=\sum a_k\lambda_k^s$ where $\lambda_1=\tdiam(C) > \lambda_2 > \cdots$ and $a_k$ is the multiplicity of $\lambda_k$, that is the number of $v\in \Vv$ with diameter $\lambda_k$.

\vspace{.3 cm}

 \subsection{The Upper Box Dimension}
 \label{cantor07.ssect-box}

\noindent This section is a reminder about the upper box dimension of a fractal.  For a treatment of the many fractal dimensions, the reader can consult \cite{FAL}.  Let $X$ be a metric space with metric $d$.  Let $N_\delta(X)$ be the least number of sets of diameter at most $\delta$ that cover $X$.

\begin{defini}
\label{cantor07.def-box}
The {\em upper box dimension} is defined as 
$$\ubox(C)=\limsup_{\delta\downarrow 0}\frac{\log N_\delta(C)}{-\log \delta}$$
\end{defini}

\noindent As shown in \cite{FAL}(Ch. 2.1), the upper box dimension satisfies the following dimension properties: monotonicity, zero on finite sets, and it gives dimension $n$ to open sets in $\RM^n$ .  Most importantly, it is invariant under bi-Lipschitz transformations.  Therefore, if two different metrics on $X$ are metrically equivalent, then they have the same upper box dimension.  The upper box dimension is also the largest of the typical fractal dimensions.  In particular, it is greater than or equal to the Hausdorff dimension of $X$.

\vspace{.3 cm}

 \subsection{The Abscissa of Convergence: Proof of Theorem~\ref{cantor07.th-box}}
 \label{cantor07.ssect-boxthm}

\noindent   In this section, the abscissa of convergence of the $\zeta$-function of $D$ will be denoted by $s_0$.  Also, $\zeta(s)$ will be written as $\sum a_k \lambda_k^s$.  In order to prove the theorem, the following classical lemma on Dirichlet series is necessary.

\begin{lemma}
\label{cantor07.lem-HaR}
Let $\zeta(s)=\sum a_k\lambda_k^s$ be a Dirichlet series with abscissa of convergence $s_0$.  Suppose further that all the $\lambda_1 > \lambda_2 > \cdots$ and that $a_k > 0$ for all $k$.  Then
$$
\limsup_{k\to\infty}\frac{\log \sum_{j=1}^{j=k} a_j}{-\log \lambda_k}=s_0
$$
\end{lemma}
\noindent {\bf Proof:}  A proof of this can be found in \cite{HaR} (Ch. 2.6).  Note that the form of the Dirichlet series used there is slightly different than the one used here.
\hfill $\Box$

\vspace{.3 cm}

\noindent With this lemma in hand, it is now possible to prove the theorem.

\vspace{.3 cm}

\noindent {\bf Proof of Theorem \ref{cantor07.th-box}:} For any $\delta > 0$ such that $\lambda_n > \delta \geq \lambda_{n+1}$, $N_\delta(C)=N_{\lambda_{n+1}}(C)$ since there are no vertices with $\delta \geq \tdiam([v]) > \lambda_{n+1}$.  Thus,

$$\frac {\log N_{\lambda_{n+1}}(C)}{-\log \lambda_{n+1}} \leq 
   \frac {\log N_{\delta}(C)}{-\log \delta} < 
    \frac {\log N_{\lambda_{n+1}}(C)}{-\log \lambda_{n}}.
$$

\noindent Let $M$ be such that every vertex has at most $M$ children.  A minimal cover of $C$ with sets of diameter at most $\lambda_n$ must use every vertex of diameter $\lambda_n$.  Thus, a cover of $C$ with sets of diameter at most $\lambda_{n+1}$ can be obtained by taking the children of each set of diameter $\lambda_n$.  This cover, $\Oo$, is in fact minimal since no $O\in\Oo$ can cover two children of a vertex of diameter $\lambda_n$.  Since every vertex of diameter at least $\lambda_n$ has at least $2$ children and at most $M_n$ children, this gives
$$
N_{\lambda_{n}} + a_n \leq N_{\lambda_{n+1}}\leq N_{\lambda_{n}} + (M_n-1)a_n.
$$
After iterating the procedure,
$$
1+\sum_{k=1}^n a_k \leq  N_{\lambda_{n+1}} \leq1+ (M_n-1) \sum_{k=1}^n a_k.
$$
where the $1$ comes from the fact that $N_{\lambda_1}=1$.  For the binary tree, it is easy to check that these inequalities are in face equalities and therefore that this estimate is in some sense optimal.  Since every cover of $C$ with sets of diameter at most $\lambda_{n+1}$ must use every vertex of diameter $\lambda_{n+1}$, then $N_{\lambda_{n+1}}\geq a_{n+1}$.  Consequently, $N_{\lambda_{n+1}}\geq 1/2(a_{n+1}+1+\sum_{k=1}^n a_k)$.  Thus, 
$$
\frac{\log 1/2(\sum_{k=1}^{j=n+1} a_k)}{-\log \lambda_{n+1}}\leq \frac{\log N_\delta(C)}{-\log \delta}< \frac{\log (1+(M_n-1)\sum_{k=1}^{j=n} a_k)}{-\log \lambda_n}
$$
Therefore, since $(\log (M_n-1))/(-\log \lambda_n)\to 0$ as $n\to \infty$ then
$$
\limsup_{n\to\infty}\frac{\log \sum_{k=1}^{j=n+1} a_k}{-\log \lambda_{n+1}}\leq \limsup_{\delta\to 0}\frac{\log N_\delta(C)}{-\log \delta}\leq \limsup_{n\to\infty}\frac{\log \sum_{k=1}^{j=n} a_k}{-\log \lambda_n}
$$
and $\ubox(C)=s_0$.
\hfill $\Box$

\vspace{.5cm}
\section{Measure Theory on $C$}
\label{cantor07.sect-meastheory}

\noindent This section extends the study of the noncommutative geometry of a Cantor set $C$ by studying a measure $\mu$ that is naturally defined on $C$.

\vspace{.3 cm}
\subsection{$\zeta$-regularity: Proof of Theorem \ref{cantor07.th-meas}}
\label{cantor07.ssect-zetareg}

\noindent  In order to study more deeply the geometry of $C$ it is necessary to make some assumptions on $C$.

\begin{defini}
\label{cantor07.def-zetareg}
A Cantor set $C$ with regular ultrametric $d$ is {\em $\zeta$-regular} if the abscissa of convergence, $s_0$, of its $\zeta$-function is finite and if for any $f\in\Cc(C)$ and any $\tau\in\Upsilon(C)$ 
\begin{equation}\label{cantor07.eqn-zetareg}
\lim_{s\downarrow s_0} (s-s_0)\TR\left( |D|^{-s} \pi_\tau(f)\right)
\end{equation}

\noindent exists.
\end{defini}

\noindent Given a $\zeta$-regular Cantor set and a choice function $\tau\in\Upsilon(C)$, it is then possible to define a measure $\mu_\tau$ on $C$ given by
$$\mu_\tau(f) \;=\; \int_C fd\mu_\tau =
   \lim_{s\downarrow s_0}
    \frac{\TR\left( |D|^{-s} \pi_\tau(f)\right)}
     {\TR\left( |D|^{-s} \right)}$$

\noindent {\bf Proof of Theorem \ref{cantor07.th-meas}:}
Let $\tau,\tau'\in\Upsilon(C)$ and $f\in \lip(C)$ with Lipschitz constant $k$.  For $\Re(s) > s_0$, since $|D|^{-s}$ and $\pi_\tau(f)$ is bounded, then $|D|^{-s}\pi_\tau(f)$ is trace class and similarly $|D|^{-s}\pi_{\tau'}(f)$ is trace class.  Therefore,
\begin{eqnarray*}
|\TR(|D|^{-s}(\pi_\tau(f)-\pi_{\tau'}(f)))|&\leq&\sum_{v\in\Vv}{|f(\tau_+(v))-f(\tau'_+(v))|}{\tdiam(v)^{\Re(s)}}\\
&&+\sum_{v\in\Vv}{|f(\tau_-(v))-f(\tau'_-(v))|}{\tdiam(v)^{\Re(s)}}\\
&\leq& 2\sum_{v\in\Vv} k \tdiam(v)^{\Re(s)+1}.
\end{eqnarray*}
Consequently,
$$|\mu_\tau(f)-\mu_{\tau'}(f)|=|\lim_{s\downarrow s_0}
    \frac{\TR\left( |D|^{-s} \pi_\tau(f)\right)-\TR\left( |D|^{-s} \pi_{\tau'}(f)\right)}
     {\TR\left( |D|^{-s} \right)}|=0$$ since $\TR(|D|^{-s_0-1})<\infty$.  Since $\lip(C)$ is dense in $\Cc(C)$ and $\pi_\tau$ is continuous for all $\tau\in\Upsilon(C)$, then $\mu_\tau$ and $\mu_{\tau'}$ are equal on $\Cc(C)$.  Since $\pi_\tau$ is faithful for all $\tau\in\Upsilon(C)$, then $\mu_\tau$ is a probability measure for each $\tau$.
\hfill $\Box$

\vspace{.3 cm}

\subsection{The Measure on the Space of Choices}
\label{cantor07.ssect-choicemeas}

\noindent In what follows, it will be necessary to have a measure on the spaces of choices, $\Upsilon(C)$.  Recall that $\Upsilon(C)$ was the set of all functions $\tau:\Vv\to C\times C$ such that $\tau(v)\in [v]\times[v]$ and $d(\tau_+(v),\tau_-(v))=\tdiam(v)$.  Let $\Gg \subset \Vv\times \Vv$ be defined to be the set of all brothers.  That is $(u,v)\in\Gg$ if $u$ and $v$ have the same parent and $u\neq v$.  Let $\Gg_v$ be the set of all brothers whose parent is $v$. Now,  $x,y\in [v]$ are such that $d(x,y)=\tdiam ([v])$ if and only if there is a unique pair $(w,w')\in \Gg_v$ of distinct children of $v$ such that $x\in [w]$ and $y\in [w']$. Consequently

$$\Upsilon(C) = 
   \prod_{v\in \Vv} 
    \bigsqcup_{(w,w')\in\Gg_v} 
     [w]\times[w'].
$$

\noindent Therefore, define a measure $\nu_v$ on $\Upsilon_v(C):=\bigsqcup_{(w,w')\in\Gg_v} [w]\times[w']$ by

$$\nu_v=
   \frac{\mu\times\mu}{\sum_{(w,w')\in\Gg_v} \mu([w])\mu([w'])}.
$$

\noindent This is then a probability measure on $\Upsilon_v(C)$.  Using the Kolmogorov Consistency theorem \cite{Par}(V.5), there is an extension of these measures to a probability measure $\nu$ on $\Upsilon(C)$.  This measure $\nu$ is such that $\nu((\prod_{w\neq v} \Upsilon_w(C))\times U_v)=\nu_v(U_v)$ for any $\nu_v$-measurable set $U_v$.

\vspace{.5 cm}
\section{Dirichlet Forms and the Operator $\Delta$}
\label{cantor07.sect-dirformsdelta}

\noindent  In this section, let $\Ll_\CM^2(C,d\mu)$ denote the Hilbert space completion of $\Cc(C,\CM)$ with respect to  $\langle f,g\rangle=\int_C \bar fgd\mu$ and let $\Ll^2(C,d\mu)$ denote the Hilbert space completion of $\Cc(C,\RM)$ with respect to the same inner product.  It is of interest to study Markovian semigroups of operators on $\Ll^2(C,d\mu)$.  As shown in \cite{Fuk}, the study of Markovian semigroups is equivalent to studying the Dirichlet forms on $\Ll^2(C,d\mu)$.

\vspace{.3 cm}
\subsection{Dirichlet Forms: Proof of Theorem \ref{cantor07.th-lapla}}
\label{cantor07.ssect-dirforms}

\noindent Given a real Hilbert space $\Hh$, a non-negative definite symmetric bilinear form densely defined on $\Hh$ is called a {\em symmetric form} on $\Hh$.  Let $Q$ be a symmetric form on a Hilbert space $H$.  If $\Dom (Q)$ is complete with respect to the metric given by $\langle f,g \rangle_1=\langle f,g \rangle_H+Q(f,g)$ where $\langle\cdot,\cdot\rangle_H$ is the inner product on $H$ then $Q$ is called a {\em closed} form.  Given a closed symmetric form $Q$ on $\Ll^2(C,d\mu)$, then $Q$ is called {\em Markovian} if $Q(\tilde f,\tilde f)\leq Q(f,f)$ where $\tilde f=\min(\max(0,f),1)$.  If $Q$ is not closed, the condition to be Markovian is more complicated; however, the previous condition is sufficient.  A closed symmetric Markovian form is called a {\em (symmetric) Dirichlet form}.  Given the formalism of the previous sections, it is possible to define a form $Q_s$ on $\Ll_\CM^2(C,d\mu)$ by 
$$Q_s(f,g):=\frac{1}{2}\int_{\Upsilon(C)} \TR(|D|^{-s}[D,\pi_\tau(f)]^*[D,\pi_\tau(g)])d\nu(\tau).$$
It is now necessary to specify a domain for the form.  Let $\Ee\subset \Ll^2(C,d\mu)$ be the real linear space spanned by $\{\chi_v : v\in \Vv\}$ where $\chi_v$ is the characteristic function of $[v]\subset C$. 

\begin{lemma}
\label{cantor07.lem-dirform}
$\Ee$ is dense in $\Ll^2(C,d\mu)$.
\end{lemma}

\noindent{\bf Proof:}
Let $f\in\Cc(C)$.  Since $f$ is continuous and $C$ is compact, then $f$ is uniformly continuous.  Consequently, for $\epsilon > 0$ there is a $\delta > 0$ such that if $d(x,y)<\delta$ then $|f(x)-f(y)| < \epsilon$.  Let $v_1,\dots,v_N$ be a partition of $C$ such that $\tdiam([v_i])<\delta$.  Let $\tau\in\Upsilon(C)$.  Then if $g(x):=f(\tau_+(v_j))$ where $v_j$ is the unique vertex of the partition such that $x\in[v_j]$.  Then $||f-g||_\infty < \epsilon$ and consequently $||f-g||_2 < \epsilon$.  Thus $\Ee$ is dense in $\Cc(C)$.  Since $\Cc(C)$ is dense in $\Ll^2(C,d\mu)$ then $\Ee$ is dense in $\Ll^2(C,d\mu)$.
\hfill $\Box$

\vspace{.3 cm}

\noindent Let $\Dom(Q_s)=\Ee$.

\vspace{.3 cm}

\noindent{\bf Proof of Theorem \ref{cantor07.th-lapla}:}
It is clear that $Q_s$ must be bilinear.  It is symmetric because of the trace and because $|D|^{-s}$ commutes with $[D,\pi_\tau(f)]$ for all $f\in\lip(C)$.  Now,
$$[D,\pi_\tau(f)]^\ast[D,\pi_\tau(g)]\psi(v)=
   \frac{f(\tau_+(v))-f(\tau_-(v))}{\tdiam(v)}
    \frac{g(\tau_+(v))-g(\tau_-(v))}{\tdiam(v)}
     \psi(v)
$$
and thus
$$\TR(|D|^{-s}[D,\pi_\tau(f)]^\ast[D,\pi_\tau(f)])=
   2\sum_{v\in\Vv} 
     \tdiam(v)^{s-2} (f(\tau_+(v))-f(\tau_-(v)))^2.
$$

\noindent Consequently, $Q_s$ is non-negative definite. Since $\chi_v(\tau_+(w))-\chi_v(\tau_-(w))= 0$  if $w \nsucc v$ then $[D,\pi_\tau(\chi_v)]$ is finite rank for each characteristic function $\chi_v$ with $v\in\Vv$.  Thus for $f\in \Ee$, $[D,\pi_\tau(f)]$ is finite rank and $Q_s(f,g) < \infty$ for all $g\in\Ll^2(C,d\mu)$.

\noindent Let now $(f_n)_{n\in\NM}$ be a sequence of functions in $\Ee$ such that $\lim_{n\to\infty} ||f_n||_{\Ll^2} = 0$ and $\lim_{n,m\to\infty} Q_s(f_n-f_m,f_n-f_m) = 0$.  To show that $Q_s$ is closable, it is then necessary to show that $\lim_{n\to\infty} Q_s(f_n,f_n) = 0$. Since $\lim_{n\to\infty} ||f_n||_{\Ll^2} = 0$ there is a subsequence $f_{n_i}$ that converges pointwise $\mu$-a.e. to $0$ \cite{Rud} (Thm. 3.12). In particular, thanks to the definition of the measure $\nu$ on the set of choices, $f_{n_i}(\tau_+(v))\to 0$ for $\nu$-a.e. choice and for all $v\in \Vv$.  Similarly for $\tau_-(v)$. So, given $\epsilon > 0$ let $N$ be such that $Q_s(f_n-f_m,f_n-f_m) < \epsilon$ for $n,m > N$.  Then for $m > N$, 
$$Q_s(f_m,f_m) = \int_{\Upsilon(C)}\sum_{j=1}^K \tdiam(v_j)^{s-2} (f_m(\tau_+(v_j))-f_m(\tau_-(v_j)))^2d\nu.$$
Since $(f_m(\tau_+(v_j))-f_m(\tau_-(v_j)))^2=$
$$\liminf_{i\to\infty} (f_m(\tau_+(v_j))-f_{n_i}(\tau_+(v_j))-f_m(\tau_-(v_j))+f_{n_i}(\tau_-(v_j)))^2$$
then using Fatou's lemma, 
$$Q_s(f_m,f_m)\leq \liminf_{i\to\infty}Q_s(f_m-f_{n_i},f_m-f_{n_i})< \epsilon.$$

\noindent Thus $\lim_{m\to\infty} Q_s(f_m,f_m)=0$ and $Q_s$ is closable.

\noindent The proof that $Q_s$ is Markovian is by inspection: let $C_-,C_0,C_+$ denote the closed subsets of $C$ for which $f\leq 0, 0\leq f\leq 1, 1\leq f$. If $\tau_+(v)\in C_i$ and $\tau_-(v)\in C_j$  then $|\tilde f(\tau_+(v))-\tilde f(\tau_-(v))|\leq |f(\tau_+(v))-f(\tau_-(v))|$ for each $i,j$.  Thus $Q_s(\tilde f,\tilde f)\leq Q_s(f,f)$.
\hfill$\Box$

\vspace{.3 cm}
 
\noindent It is now possible to get a closed Dirichlet form using the following result.

\begin{theo}[\cite{Fuk} Thm 2.1.1]
\label{cantor07.th-dirform2}
Suppose Q is a closable Markovian symmetric form on $\Ll^2(X,m)$ where $X$ is a locally compact separable Hausdorff space and $m$ is a positive Radon measure on $X$ such that $\Supp(m)=X$.  Then its smallest closed extension is a Dirichlet form. 
\end{theo}
 
\vspace{.3cm}

\subsection{Self-Adjoint Operators and Operator Semigroups}
\label{cantor07.ssect-opersemigroup}

\noindent  This section follows \cite{Fuk} (Ch 1.3).  Let $H$ be a real Hilbert space.

\begin{defini}
\label{cantor07.def-semigroup}
A family $\{T_t, t>0\}$  of linear operators is called a {\em strongly continuous, symmetric, contraction semigroup} if:

\noindent (i) each $T_t$ is a symmetric operator with $\Dom(T_t)=H$.

\noindent (ii) semigroup property: $T_tT_s=T_{t+s}$ for $t,s>0$.

\noindent (iii) contraction property: $\langle T_tf,T_tf\rangle \leq \langle f,f\rangle$ for all $f\in H$ and $t > 0$.

\noindent (iv) strong continuity: $\langle T_tf-f,T_tf-f\rangle \to 0$ as $t\downarrow 0$ for all $f\in H$.
\end{defini}

\noindent Let $\{T_t, t>0\}$ be such a semigroup.  Then the {\em generator} $A$ is an operator on $H$ defined by
$$Af:=\lim_{t\downarrow 0} \frac{T_tf-f}{t}, \Dom(A):= \{f\in H: Af \text{ exists as a strong limit}\}.$$
In fact, there is a one-to-one correspondence between non-positive definite self-adjoint operators on $H$ and the family of strongly continuous, symmetric, contraction semigroups.  The correspondence from $A$ to $\{T_t\}$ is given by $T_t=\exp(tA)$.

\noindent Given a non-positive definite self-adjoint operator, let $Q(u,v):=\langle -Au,u \rangle$ with $\Dom(Q):=\Dom(\sqrt{-A})$.  It turns out that $Q$ is a closed symmetric form on $H$.  This correspondence is also one-to-one.  Starting with a closed, symmetric form $Q$ on $H$ the construction of $A$ is slightly more involved.    Since $Q$ is closed, then $\Dom(Q)$ is a Hilbert space with norm $||g||_1=||g||_{\Ll^2}+Q(g,g)$.  Fix $f\in H$.  Then $\langle \cdot, f\rangle$ is a bounded linear functional on $\Dom(Q)$.  Therefore, let $Bf$ be the unique vector in $\Dom(Q)$ corresponding to this linear functional by the Riesz Representation Theorem.  Let $A:= I- B^{-1}$.  Then $A$ is the non-positive definite self-adjoint operator corresponding to $Q$.

\noindent Now let $H=\Ll^2(X,m)$ where $X$ is a locally compact separable Hausdorff space and $m$ is a positive Radon measure on $X$ such that $\Supp(m)=X$.  A bounded linear operator $S$ on $\Ll^2(X,m)$ is called {\em Markovian} if $0\leq Sf \leq 1, m$-a.e. whenever $f\in \Ll^2(X,m)$ is such that $0\leq f\leq 1$. A strongly continuous, symmetric, contraction semigroup $\{T_t\}$ such that $T_t$ is Markovian for each $t > 0$ is called a {\em Markovian semigroup}.

\begin{theo}[\cite{Fuk} Thm 1.4.1]
\label{cantor07.th-dirformsemi}
Let $X$ be a locally compact separable Hausdorff space and $m$ a positive Radon measure on $X$ such that $\Supp(m)=X$.  Then there is a one-to-one correspondence between Dirichlet forms on $\Ll^2(X,m)$ and Markovian semigroups on $\Ll^2(X,m)$.
\end{theo}

\vspace{.3cm}

\subsection{The Operators $\Delta_s$}
\label{cantor07.ssect-delta}

\noindent  Let $C$ be a $\zeta$-regular Cantor set with regular ultrametric $d$.  Let $\mu$ be the measure constructed via the $\zeta$-function.  Suppose $\mu$ is such that $\Supp(\mu)=C$.  Then for $s\in\RM$, the previous results give a non-positive definite self-adjoint operator $\Delta_s$ such that $T_t:=\exp(t\Delta_s)$ is a Markovian semigroup.  $\Delta_s$ is such that

$$\langle -\Delta_s f,g\rangle = 
   \frac{1}{2}\int_{\Upsilon(C)} 
    \TR(|D|^{-s}[D,\pi_\tau(f)]^\ast[D,\pi_\tau(g)])
      d\nu(\tau)
$$

\noindent for $f,g\in\Dom(\Delta_s)$.  It is important to note that $\Ee \subset \Dom(\Delta_s)\subset \Dom(\bar Q_s)$ where $\bar Q_s$ is the smallest closed extension of $Q_s$.  

\noindent It is possible to calculate $\Delta_s\chi_v$ for $v\in \Vv$.  Since $\chi_v(\tau_+(w))-\chi_v(\tau_-(w))= 0$  if $w \nsucc v$, then for $g\in\Dom(Q_s)$, $ \langle -\Delta_s \chi_v, g \rangle$ = 
$$\sum_{w\succ v} \tdiam(w)^{s-2} \int_{\Upsilon(C)} (\chi_v(\tau_+(w))-\chi_v(\tau_-(w)))(g(\tau_+(w))-g(\tau_-(w)))d\nu(\tau).$$
Since $\tau$ is only applied to $w$, then by the very definition of $\nu_w$
$$=\sum_{w\succ v} \frac{\tdiam(w)^{s-2}}{\sum_{(u,u')\in \Gg_w} \mu([u])\mu([u'])} \sum_{(u,u')\in \Gg_w}\int_{[u]\times[u']} (\chi_v(x)-\chi_v(y))(g(x)-g(y))d\mu d\mu.$$
For $w$ an ancestor of $v$ let $u_v$ be its child that is also an ancestor of $v$.  Then for any other child $u$ of $w$, $\chi_u(x)=0$ for $x\in[u]$.  Thus since $\bigcup_{(u_v,u')\in\Gg_w}[u']=[w]\cap[u_v]^c$, then 
$$=\sum_{w\succ v} \frac{\tdiam(w)^{s-2}}{\sum_{(u,u')\in \Gg_w} \mu([u])\mu([u'])} 2\int_{[v]}d\mu(x)\int_{[w]\cap[u_v]^c}g(x)-g(y)d\mu(y).$$
Consequently,
\begin{equation}
\label{cantor07.eq-deltachi}
\Delta_s \chi_v=-\sum_{w\succ v} \frac{\tdiam(w)^{s-2}}{\sum_{(u,u')\in \Gg_w} \mu([u])\mu([u'])} 2(\mu([w]\cap[u_v]^c)\chi_v-\mu([v]) \chi_{[w]\cap[u_v]^c}).
\end{equation}

\noindent An application of this formula is given by the following:

\begin{proposi}
\label{cantor07.prop-purepoint}
The spectrum of $\Delta_s$ is pure point.
\end{proposi}

\noindent{\bf Proof:}
Let $\Ll_n\subset \Ll^2(C,d\mu)$ be the space spanned by all $\chi_v$ such that $\theight(v)\leq n$.  Since $\Tt$ is Cantorian then $\dim(\Ll_n)<\infty$.  Moreover, $\Ll_n\subset \Ll_{n+1}$ and $\bigcup_n \Ll_n$ is dense in $\Ll^2(C,d\mu)$.  Equation \ref{cantor07.eq-deltachi} then gives that $\Delta_s$ leaves each $\Ll_n$ invariant.  Since $\Delta_s$ restricted to each finite dimensional $\Ll_n$ is pure point, then $\Delta_s$ is pure point.
\hfill $\Box$

\vspace{.5 cm}
\section{The Triadic Cantor Set}
\label{cantor07.sect-triadic}

\vspace{.3 cm}
\subsection{Eigenvalues and Eigenstates for $\Delta_s$ on $C_3$}
\label{cantor07.ssect-eig}

\noindent  This section will apply much of the previous machinery to the triadic Cantor set.  Let $C_3$ denote the triadic Cantor set seen as a subset of the interval $[0,1]$.  As seen in Example \ref{cantor07.exam-bin}, $C_3$ is the boundary of the infinite binary tree $\pT_2$ and has a natural homeomorphism with $\{0,1\}^\NM$ by

$$\phi(\omega)=\sum_{n=0}^\infty \frac{2\omega_n}{3^{n+1}},\ \ \ \ \omega=\{\omega_n\}_{n\in\NM}\in \{0,1\}^\NM.
$$

\noindent Let $d$ be the regular ultrametric corresponding to the weight $\epsilon(v)=3^{-\theight(v)}$.  Then for $x,y\in C_3$,

$$\frac{d(x,y)}{3}\leq |x-y|\leq d(x,y)
$$

\noindent and thus $d$ is metrically equivalent to the Euclidean metric.  Then 

$$\zeta(s)=\sum_{n=0}^\infty (\frac{2}{3^s})^n
$$

\noindent and therefore has abscissa of convergence $s_0=\ln 2/\ln 3$.  This pole is clearly a simple pole. For any $v\in\Vv$,

$$\frac{1}{2}\TR(|D|^{-s}\pi_\tau(\chi_v))=\sum_{w\preceq v} \tdiam(w)^s=\tdiam(v)^s\zeta(s)
$$

\noindent since the subtree starting at $v$ is identical to the tree starting at the root.  Consequently, $\mu(\chi_v)=\tdiam(v)^{s_0}$. Thus $\mu(f)$ is defined on all characteristic functions and can be extended to all continuous functions. Therefore, $C_3$ is $\zeta$-regular and

$$\mu([v])=\tdiam(v)^{s_0}=\frac{1}{3^{s_0\mheight(v)}}=\frac{1}{2^{\mheight(v)}}.
$$

\noindent Since $\Supp(\mu)=C_3$ then $\Delta_s$ can be defined on $\Ll^2(C_3,d\mu)$. Equation \ref{cantor07.eq-deltachi} then gives that for $v=v_0\cdots v_n\in\Vv$ with $n\geq 1$,

$$ \Delta_s \chi_v=-\sum_{j=0}^{n-1} \frac{3^{j(2-s)}}{2^{-(2j+1)}} 2(2^{-(j+1)}\chi_v-2^{-n} \chi_{[w]\cap[u_v]^c}).
$$

\noindent Letting $\bar a=1-a$ for $a\in\{0,1\}$ then this becomes
\begin{equation}\label{cantor07.eq-cantdelt}
\Delta_s \chi_v=-2\sum_{j=0}^{n-1} \left(\frac{2}{3^{s-2}}\right)^j\chi_v+ \frac{4}{2^n}\sum_{j=0}^{n-1} \left(\frac{4}{3^{s-2}}\right)^j\chi_{v_0\cdots v_{j}\bar v_{j+1}}.
\end{equation}
This formula can be used to find the eigenstates of $\Delta_s$.
\begin{defini}
\label{cantor07.def-haar}
Let $\Ww$ be the set of infinite sequences $\omega=\omega_1\omega_2\cdots\in\{0,1\}^{\NM^+}$ such that all but a finite number of $\omega_k$'s are $0$.  Let $|\omega|$ be the maximum integer $k$ such that $\omega_k=1$ with the convention that $|\omega|=0$ if $\omega=00\cdots$.  The {\em Haar function} $\phi_\omega$ is defined by
$$\phi_\omega=\sum_{v\in\{0,1\}^n} (-1)^{\omega\cdot v}\chi_v,\hspace{.5 cm} \omega\cdot v=\sum_{k=1}^n\omega_kv_k.$$
for any $n\geq |\omega|$.
\end{defini}

\noindent    Because $\chi_{v_1\dots v_N 0}+\chi_{v_1\dots v_N1}=\chi_{v_1\dots v_N}$ and since if $N=|\omega|$ then $\omega_{N+m}=0$ for $m > 0$,then
$$\sum_{v\in\{0,1\}^{N+1}} (-1)^{\omega\cdot v}\chi_v=\sum_{v\in\{0,1\}^{N}} (-1)^{\omega\cdot v}(\chi_{v_1\dots v_N0}+\chi_{v_1\dots v_N1})=\sum_{v\in\{0,1\}^n} (-1)^{\omega\cdot v}\chi_v.$$
Therefore, $\phi_\omega$ does not depend on the choice of $n$ and $\phi_\omega$ is well-defined.  Moreover, it is straightforward to check that the Haar functions are orthonormal in the sense that $\langle \phi_\omega,\phi_\sigma \rangle = \delta_{\omega,\sigma}$
for $\omega,\sigma\in\Ww$.  In addition, 
$$\chi_v=\frac{1}{2^n}\sum_{u\in\{0,1\}^n} (-1)^{v\cdot u}\phi_{u00\cdots}$$
 for $v\in\Vv$ and thus the Haar functions are an orthonormal basis for $\Ll^2(C_3,d\mu)$.  The importance of the Haar functions comes from the following theorem.

\begin{theo}
\label{cantor07.th-haar}
Let $C_3$ be the triadic Cantor set with the regular ultrametric $d$ given above.  Let $\mu$ be its associated measure.  Then

\noindent (i) The eigenstates of $\Delta_s$ are given by the Haar functions $\phi_\omega$ with $\omega\in\Ww$.

\noindent (ii) The eigenvalues of $\Delta_s$ are given by $\lambda_0=0$ and for $n\geq 1$
$$-\lambda_n=-2\left(1+3^{s_0+2-s}+\cdots+\left(3^{s_0+2-s}\right)^{n-2}+2\left(3^{s_0+2-s}\right)^{n-1}\right)$$

\noindent (iii) The degeneracy of $\lambda_n$ is $2^{n-1}$ for $n\geq 1$ whereas $\lambda_0$ is simple.

\noindent (iv) For $s > s_0+2$, $\Delta_s$ is bounded and is a compact perturbation of a multiple of the identity.

\noindent (v) For $s\leq s_0+2$, $\Delta_s$ has compact resolvent.

\noindent (vi) For $s < s_0 +2$, the density of states $\Nn(\lambda)$ given by the dimension of the spectral space corresponding to eigenvalues whose magnitude is less than or equal to $\lambda$ satisfies
$$\Nn(\lambda)\stackrel{\lambda\uparrow \infty}{\sim} 2\left(\frac{\lambda}{2k}\right)^{{s_0}/(2+s_0-s)}(1+o(1))$$
where $k=1/(1-3^{s-2-s_0})+1$.
\end{theo}

\begin{rem}
\label{cantor07.rem-weyl}
{\em On a compact Riemannian manifold $M$, the Laplacian is an unbounded operator with compact resolvent.  Moreover, Weyl's theorem says that if $m$ is the dimension of $M$ then $\Nn(\lambda)\sim c_0 \lambda^{m/2}$ as $\lambda\to\infty$ for an appropriate constant $c_0$.  The constant $c_0$ is not arbitrary and actually gives the volume of the unit ball in the cotangent bundle over the manifold.  In any case, the previous theorem shows that if $\Delta_s$ is interpreted as the Laplacian on a compact Riemannian manifold then $m=2s_0/(2+s_0-s)$ gives the Riemannian dimension of this noncommutative manifold.  By analogy, this suggests that $\Delta_{s_0}$ is the appropriate Laplacian on $C_3$ since it gives Riemannian dimension $s_0$. }
\end{rem}

\noindent{\bf Proof:}
Using Equation \ref{cantor07.eq-cantdelt} and the definition of the Haar function, for $\omega\in\Ww$ with $|\omega|=n > 0$
\begin{eqnarray*}
-\Delta_s\phi_\omega&=&\sum_{v\in\{0,1\}^n} (-1)^{\omega\cdot v}\left(2\sum_{j=0}^{n-1} \left(\frac{2}{3^{s-2}}\right)^j\chi_v- \frac{4}{2^n}\sum_{j=0}^{n-1} \left(\frac{4}{3^{s-2}}\right)^j\chi_{v_0\cdots v_{j}\bar v_{j+1}}\right)\\
&=&2\sum_{j=0}^{n-1} \left(\frac{2}{3^{s-2}}\right)^j\phi_\omega-\frac{4}{2^n}\sum_{j=0}^{n-1}\left(\frac{4}{3^{s-2}}\right)^j\sum_{v\in\{0,1\}^n} (-1)^{\omega\cdot v} \chi_{v_0\cdots v_{j}\bar v_{j+1}}.
\end{eqnarray*}
For $j < n-1$ the last sum on the right hand side vanishes and for $j=n-1$
$$\sum_{v\in\{0,1\}^n} (-1)^{\omega\cdot v} \chi_{v_0\cdots v_{j}\bar v_{j+1}}=-\phi_\omega$$ since $(-1)^{v_n}=-(-1)^{\bar v_n}$.
Consequently, 
$$\Delta_s\phi_\omega=-\left(2\sum_{j=0}^{n-1} \left(3^{s_0+2-s}\right)^j+2\left(3^{s_0+2-s}\right)^{n-1}\right)\phi_\omega.$$
Therefore, the Haar basis is an eigenbasis for $\Delta_s$ and the corresponding eigenvalues are precisely the $-\lambda_n$'s given in the statement of the theorem.  Since there are exactly $2^{n-1}$ sequences $\omega\in \Ww$ with $|\omega|=n$ for $n > 0$ then the degeneracy of $-\lambda_n$ is $2^{n-1}$.

\noindent If $3^{s_0+2-s} < 1$, that is if $s > s_0 +2$ then as $n\to \infty$,
$$-\lambda_n=-2\sum_{j=0}^{n-1}\left(3^{s_0+2-s}\right)^j+2\left(3^{s_0+2-s}\right)^{n-1}\to-\frac{2}{1-3^{s_0+2-s}}=:-\lambda_\infty.$$
Hence, $\Delta_s$ is bounded and $\Delta_s+\lambda_\infty\id$ is compact.

\noindent If $s=s_0+2$ then $3^{s_0+2-s}=1$ and $-\lambda_n=-2(n+1)$.  Therefore, $(\Delta_s^2+1)^{-1}$ is compact and $\Delta_s$ has compact resolvent.  If $s < s_0+2$ then $3^{s_0+2-s}>1$ and
$$-\lambda_n=-2\left(3^{s_0+2-s}\right)^{n-1}\left(\frac{1-(3^{s-2-s_0})^n}{1-3^{s-2-s_0}}+1\right).$$
Therefore, $(\Delta_s^2+1)^{-1}$ is compact and $\Delta_s$ has compact resolvent.  Moreover,
 if $N(\lambda)$ is such that
$$\lambda=2\left(3^{s_0+2-s}\right)^{N(\lambda)-1}\left(\frac{1-(3^{s-2-s_0})^{N(\lambda)}}{1-3^{s-2-s_0}}+1\right)$$
then if $k:=1/(1-3^{s-2-s_0})+1$, 
$$N(\lambda)=1+\frac{\ln(\lambda+2(3^{s_0+2-s}-1)^{-1})-\ln 2k}{\ln 2-(s-2)\ln 3}.$$ 
Now,  
$$\lim_{\lambda\to\infty}(N(\lambda)- \frac{\ln (\lambda/(2k))}{\ln 2-(s-2)\ln 3})=0.$$
Thus, since 
$$\Nn(\lambda)=1+\sum_{n\geq 1,\lambda_n\leq \lambda} 2^{n-1}=2^{N(\lambda)}$$
then
$$\Nn(\lambda)\sim 2\left(\frac{\lambda}{2k}\right)^{{s_0}/(2+s_0-s)}(1+o(1))$$
as $\lambda\to\infty$ as desired.
\hfill $\Box$

\vspace{.3 cm}
\subsection{Diffusion on $C_3$}
\label{cantor07.ssect-diffusion}

\noindent Having computed the eigenstates and eigenvalues of $\Delta_s$, it is now possible to get an explicit description of its associated Markovian semigroup $\{\exp(t\Delta_s)\}_{t>0}$. In order to do so, let
$$\kappa_n(x,y):=\left\{ \begin{array}{ccc}1 & \text{if } d(x,y)=3^{-n}\\ 0 & \text{otherwise}\end{array}\right..$$

\begin{theo}
\label{cantor07.theo-heat}
Under the assumptions of Theorem \ref{cantor07.th-haar} and for $s < s_0+2$, the following hold

\noindent (i) Let the heat kernel $K_t(x,y)$ be defined by
$$\langle f,e^{t\Delta_s}g\rangle = \int_{C_3\times C_3} f(x)K_t(x,y)g(y)d\mu(x)d\mu(y)$$
for $f,g\in\Ll^2(C_3,d\mu)$.  Then, $K_t(x,y)=\sum_{n=0}^\infty \kappa_n(x,y)a_n(t,s)$
where $$a_n(t,s)=1-2^{n}e^{-t\lambda_{n+1}}+\sum_{m=1}^n2^{m-1}e^{-t\lambda_m}$$ for $n\geq 1$ and $a_0=1-e^{-t\lambda_1}$.  Moreover, $K_t(x,y)$ is positive for all $x,y\in C_3$ and $t > 0$.  In addition, $K_t\in\Ll^\infty(C_3\times C_3,\mu\times\mu)$ for $t > 0$.

\noindent (ii) The Markovian semigroup $\{e^{t\Delta_s}\}$ defines a Markov process $(X_t)_{t\geq 0}$ with values in $C_3$ defined by
$$\EM(f_1(X_{t_1})\cdots f_n(X_{t_n}))=\langle 1,\hat f_ne^{(t_n-t_{n-1})\Delta_s}\cdots\hat f_1e^{t_1\Delta_s}1\rangle$$
where $f_k\in\Cc(C)$, $\hat f$ denotes the operator on $\Ll^2(C_3,d\mu)$ given by multiplication by $f$, and where $t_n > \cdots > t_1 > 0$.  This Markov process is stationary and satisfies the following for $s$ fixed:
$$\EM(d(X_{t_0},X_{t_0+t})^\beta)\stackrel{t\downarrow 0}{\sim}$$$$ \left(\frac{\lambda_1}{2}+\frac{1}{2}\sum_{n=1}^\infty \left(\frac{1}{2\cdot 3^\beta}\right)^n\left(2^{n}\lambda_{n+1}-\sum_{m=1}^n2^{m-1}\lambda_m\right)\right)t(1+o(1))$$
for $\beta > s_0+2-s$ and
$$\EM(d(X_{t_0},X_{t_0+t})^\beta)\stackrel{t\downarrow 0}{\sim}\frac{1}{2\beta\ln 3}\left(\frac{1}{1-3^{-\beta}}+1\right)\left(1-\frac{1}{3^{\beta+s_0}-1}\right) t\ln(1/t)\left(1+o(1)\right)$$
for $\beta=s_0+2-s$.  For $\beta < s_0+2-s$, $$\EM(d(X_{t_0},X_{t_0+t})^\beta)=O(t^{\beta/(s_0+2-s)}\ln(1/t)).$$ 
\end{theo}

\begin{rem}
\label{cantor07.rem-diffusion}
{\em The previous section had suggested that $\Delta_{s_0}$ is the proper generalization of the Laplacian to the Cantor set.  Classical Brownian motion on the real line is generated by the Laplacian and satisfies $\EM(|X_{t_0}-X_{t+t_0}|^2)=|t|$.  For $s=s_0$, $\EM(d(X_{t_0},X_{t_0+t})^2){\sim}t\ln(1/t)$ and so there is a subdominant contribution by a term of order $\ln(1/t)$.  For $\beta=2$ this subdominant contribution only appears for $s\leq s_0$ and therefore suggests that on the Cantor set something special is happening at $s=s_0$ as the subdominant term $t\ln(1/t)$ takes over from the term $t$ which dominates for $s > s_0$.  A further understanding of this phenomenon needs to be investigated although presumably this logarithmic singularity comes from the fact that $X_t$ describes a jump process across the gaps of the Cantor set. }
\end{rem}

\noindent {\bf Proof:}
Because of the spectral decomposition of $\Delta_s$ given in Theorem \ref{cantor07.th-haar}, 
$$e^{t\Delta_s}=\sum_{n=0}^\infty e^{-t\lambda_n}\Pi_n$$
where $\Delta_s$ is the spectral projection onto the eigenspace of $\Delta_s$ corresponding to the eigenvalue $-\lambda_n$.  For $n=0$, $\Pi_0=|\phi_{00\cdots}\rangle\langle\phi_{00\cdots}|=|1\rangle\langle 1|$.  For $n\geq 1$,
$$\Pi_n=\sum_{\sigma_1,\dots,\sigma_{n-1}\in\{0,1\}} |\phi_{\sigma_1\cdots\sigma_{n-1}100\cdots}\rangle \langle\phi_{\sigma_1\cdots\sigma_{n-1}100\cdots}|$$
since $\phi_{\sigma_1\cdots\sigma_{n-1}100\cdots}$ generate the eigenspace of Haar functions $\phi_\omega$ with $|\omega|=n$.  By the definition of the Haar function,
\begin{eqnarray*}
\Pi_n&=&\sum_{u_k,v_k\in\{0,1\},k=1,\dots,n}(-1)^{u_n-v_n}|\chi_u\rangle\langle \chi_v|\sum_{\sigma_1,\dots,\sigma_{n-1}\in\{0,1\}}\prod_{k=1}^{n-1}(-1)^{(u_k-v_k)\sigma_k}\\
&=&2^{n-1}\sum_{u\in\{0,1\}^{n-1}} |\chi_{u0}\rangle\langle\chi_{u0}|-|\chi_{u0}\rangle\langle\chi_{u1}| -|\chi_{u1}\rangle\langle\chi_{u0}|+|\chi_{u1}\rangle\langle\chi_{u1}|.
\end{eqnarray*}
Now $|\chi_{u}\rangle\langle\chi_{v}|$ is the operator with functional kernel $\chi_u(x)\chi_v(y)$.  Because
$$\sum_{u\in\{0,1\}^{n-1}} \chi_{u0}(x)\chi_{u0}(y)+\chi_{u1}(x)\chi_{u1}(y)=\left\{ \begin{array}{cl}1 & \text{if } d(x,y)\leq 3^{-n}\\ 0 &\text{otherwise}\end{array}\right.$$
and
$$\sum_{u\in\{0,1\}^{n-1}} \chi_{u0}(x)\chi_{u1}(y)+\chi_{u1}(x)\chi_{u0}(y)=\left\{ \begin{array}{cl}1 & \text{if } d(x,y)= 3^{-n+1}\\ 0 & \text{otherwise}\end{array}\right.$$
then 
$$\Pi_n(x,y)=\left\{ \begin{array}{cl}2^{n-1}& \text{if } d(x,y)\leq3^{-n}\\-2^{n-1} & \text{if } d(x,y)= 3^{-n+1}\\ 0 & \text{otherwise}\end{array}\right.$$
where $\Pi_n(x,y)$ is the functional kernel of the operator $\Pi_n$.  Using the functions $\kappa_n$, this becomes
$$\Pi_n(x,y)=2^{n-1}(-\kappa_{n-1}(x,y)+\sum_{m\geq n} \kappa_m(x,y)).$$
Therefore,
\begin{eqnarray*}
K_t(x,y)&=&\sum_{n=0}^\infty e^{-t\lambda_n}\Pi_n(x,y)\\
&=&\kappa_0(x,y)(1-e^{-t\lambda_1})\\
& &+\sum_{n=1}^\infty \kappa_n(x,y)\left(1-2^{n}e^{-t\lambda_{n+1}}+\sum_{m=1}^n2^{m-1}e^{-t\lambda_m}\right).
\end{eqnarray*}
The convergence of $K_t(x,y)$ in $\Ll^2(C_3\times C_3,\mu\times \mu)$ is shown as follows.  To begin,
$$\kappa_n(x,y)=\sum_{v\in\{0,1\}^n} \chi_{v0}(x)\chi_{v1}(y)+ \chi_{v1}(x)\chi_{v0}(y)$$
gives that 
$$\int_{C_3\times C_3} \kappa_n(x,y)^2d\mu(x)d\mu(y)=\sum_{v\in\{0,1\}^n}\frac{2}{2^{2n+2}}=\frac{1}{2^{n+1}}.$$
Therefore, the corresponding norm in $\Ll^2(C_3\times C_3,\mu\times \mu)$ is $||\kappa_n||_2=2^{-(n+1)/2}$.  The coefficients of the $\kappa_n$'s in $K_t$ are positive for $t > 0$ since 
$$a_n(t,s)=2^{n}(1-e^{-t\lambda_{n+1}})-\sum_{m=1}^n2^{m-1}(1-e^{-t\lambda_m}) >1- e^{-t\lambda_{n+1}}> 0.$$
It is straightforward to show that for $t > 0$,
$$1 + \sum_{m=1}^\infty2^{m-1}e^{-t\lambda_m} =M_s(t)< \infty$$
for all $s < s_0+2$ and thus $0 < a_n(t,s)  < M_s(t).$ Since the $\kappa_n$'s have disjoint support, then $K_t$ is bounded and $K_t\in\Ll^\infty(C_3\times C_3,\mu\times\mu)$.

The definition of the stochastic process $(X_t)_{t\geq 0}$ is standard and results from the Chappman - Kolmogorov Equations.  It gives a way to evaluate $\EM(d(X_{t_0},X_{t_0+t})^\beta)$ by
$$\EM(d(X_{t_0},X_{t_0+t})^\beta)=\int_{C_3\times C_3}K_t(x,y)d(x,y)^\beta d\mu(x)d\mu(y).$$
Thus
$$\EM(d(X_{t_0},X_{t_0+t})^\beta)=\frac{1}{2}\sum_{n=0}^\infty \left(\frac{1}{2\cdot 3^\beta}\right)^na_n(t,s).$$
Now for $t > 0$ and $\beta > s_0+2-s$,
$$\frac{1}{t}\EM(d(X_{t_0},X_{t_0+t})^\beta)\leq \frac{1}{2t}\sum_{n=0}^\infty \frac{1}{3^{\beta n}} (1-e^{-t\lambda_{n+1}})\leq \frac{1}{2}\sum_{n=0}^\infty \frac{1}{3^{\beta n}} \lambda_{n+1} < \infty$$
and therefore by dominated convergence,
$$\lim_{t\to 0}\frac{\EM(d(X_{t_0},X_{t_0+t})^\beta)}{t}=\frac{\lambda_1}{2}+\frac{1}{2}\sum_{n=1}^\infty \left(\frac{1}{2\cdot 3^\beta}\right)^n\left(2^{n}\lambda_{n+1}-\sum_{m=1}^n2^{m-1}\lambda_m\right)$$
and this limit is positive and finite.  For $\beta = s_0+2-s$, let $N_t=\ln(1/t)/(\beta\ln 3)$.  First of all,
$$\frac{1}{2}\sum_{n=N_t+1}^\infty \left(\frac{1}{2\cdot 3^\beta}\right)^na_n(t,s)<\frac{1}{2}\sum_{n=N_t+1}^\infty \left(\frac{1}{3^\beta}\right)^n=\frac{t}{2\cdot 3^\beta}\frac{1}{1-3^{-\beta}}$$
and 
\begin{eqnarray*}
\frac{1}{2}\sum_{n=N_t+1}^\infty \left(\frac{1}{2\cdot 3^\beta}\right)^na_n(t,s)&>&\frac{1}{2}(1-e^{-t\lambda_{N_t+2}})\sum_{n=N_t+1}^\infty \left(\frac{1}{2\cdot3^\beta}\right)^n\\
&>&(1-e^{-t\lambda_{1}})\frac{t^{1+s_0/\beta}}{4\cdot 3^\beta}\frac{1}{1-2^{-1}3^{-\beta}}.
\end{eqnarray*}
By taking a Taylor expansion,
$$2^nt\lambda_{n+1}-2^n\frac{t^2\lambda_{n+1}^2}{2}-t\sum_{m=1}^n2^{m-1}\lambda_m\leq a_n(t,s)$$
and
$$a_n(t,s)\leq 2^nt\lambda_{n+1}-t\sum_{m=1}^n2^{m-1}\lambda_m-\frac{t^2}{2}\sum_{m=1}^n2^{m-1}\lambda_m^2.$$
Now
$$\frac{1}{2}\sum_{n=1}^{N_t} \left(\frac{1}{2\cdot3^\beta}\right)^n 2^n t^2\lambda_{n+1}^2 <t^2(\frac{1}{1-3^{-\beta}}+1)\sum_{n=1}^{N_t} 3^{n\beta}< c_0t^23^{\beta N_t}=c_0t $$
for some constant $c_0 > 0$.  Similarly, there exists $c_1 > 0$ such that
$$\sum_{n=1}^{N_t}\left(\frac{1}{2\cdot3^\beta}\right)^n\frac{t^2}{2}\sum_{m=1}^n2^{m-1}\lambda_m^2 < c_1 t.$$
Since,
$$\left(\frac{1}{2\cdot3^\beta}\right)^n\left(2^n\lambda_{n+1}-\sum_{m=1}^n2^{m-1}\lambda_m\right)=$$$$\left(\frac{1}{1-3^{-\beta}}+1\right)\left(1-\sum_{m=1}^n(2\cdot 3^\beta)^{-m}\right)-\frac{3^{-\beta}}{1-3^{-\beta}}\left(\frac{1}{2\cdot3^\beta}\right)^n$$
then
$$\frac{t}{2}\sum_{n=1}^{N_t}\left(\frac{1}{2\cdot3^\beta}\right)^n\left(2^n\lambda_{n+1}-\sum_{m=1}^n2^{m-1}\lambda_m\right)=$$$$\frac{1}{2}\left(\frac{1}{1-3^{-\beta}}+1\right)\left(1-\frac{1}{3^{\beta+s_0}-1}\right)tN_t+c_2t(1-t^{1+s_0/\beta})$$
where $c_2 > 0$ is a constant.
Consequently,
$$\EM(d(X_{t_0},X_{t_0+t})^\beta)\stackrel{t\downarrow 0}{\sim}\frac{1}{2\beta\ln 3}\left(\frac{1}{1-3^{-\beta}}+1\right)\left(1-\frac{1}{3^{\beta+s_0}-1}\right) t\ln(1/t)\left(1+o(1)\right)$$
for $\beta=s_0+2-s$.  The proof that $\EM(d(X_{t_0},X_{t_0+t})^\beta)=O(t^{\beta/(s_0+2-s)}\ln(1/t))$ is the same as above using the fact that $1-e^{-x}\leq x^{\alpha}$ for $0\leq \alpha\leq 1$.
\hfill$\Box$

\vspace{.3 cm}

\subsection{Relationship with the Vladimirov Operator}
\label{cantor07.ssect-vladimirov}

Let $p$ be a prime number.  It is a basic fact (see \cite{Sch}) that the Cantor set is homeomorphic to the $p$-adic integers, $\ZM_p$.  In fact, $\ZM_p$ is the boundary of the tree $\Tt_p$ where every vertex has exactly $p$ children and the weight function is $\epsilon(v)=p^{-\mheight(v)}$ for $v\in\Vv$.  The $p$-adic numbers are the completion of $\QM$ with respect to this ultrametric $|\cdot|_p$ and $\ZM_p$ is then the closed unit disc in $\QM_p$.  The Vladimirov operator \cite{VVZ} is constructed using the field structure of $\QM_p$.  It is defined by 
$$(\Dd \psi)(x)=\frac{p^2}{p+1}\int_{\QM_p} \frac{\psi(x)-\psi(y)}{|x-y|_p^2}dy$$
where $\psi:\QM_p\to \RM$ is a locally constant function with compact support and the measure $dy$ is the Haar measure on $\QM_p$. 

\begin{proposi}
\label{cantor07.prop-vladimirov}
For $z=v_0v_1\cdots\in\pT_2$ and $f\in\Ee$, $$(\Dd f)(z)=\frac{1}{3}\lim_{n\to\infty} \frac{1}{\mu([v_n])}\langle \chi_{v_n},-\Delta_2 f\rangle.$$
\end{proposi}

\noindent {\bf Proof:}
From Section \ref{cantor07.ssect-delta}, 
$$\langle \chi_{v_n},-\Delta_2 f\rangle = \int_{[v_n]}d\mu(x) \sum_{j=0}^{n-1} \frac{4}{\mu([v_0\cdots v_{j-1}])^2}\int_{[v_0\cdots v_{j-1}\bar v_j]}f(x)-f(y)d\mu(y).$$
But for $x\in[v_0\cdots v_j]$ and $y\in[v_0\cdots v_{j-1}\bar v_j]$, $|x-y|_2=\mu([v_0\cdots v_{j-1}])$.  Therefore,
$$\langle \chi_{v_n},-\Delta_2 f\rangle=4\int_{[v_n]}d\mu(x)\int_{[v_n]^c}\frac{f(x)-f(y)}{|x-y|_2^2}d\mu(y)$$
and the result follows.
\hfill $\Box$

\vspace{.3 cm}

\noindent Because $|D|^{-1}D=F$ is the phase of the operator $D$, then this result shows that since the Vladimirov operator is constructed out of the phase then it does not take the metric on $C_3$ into account.  This makes sense because the Vladimirov operator was created using the $2$-adic metric which comes from the measure and not from the metric on $C_3$.
\vspace{.5 cm}

\section{Conclusion and Open Problems}
\label{cantor07.sect-conclusion}

\noindent The present paper has constructed the appropriate machinery from Noncommutative Geometry to investigate various examples of ultrametric Cantor sets as Noncommutative Riemannian spaces.  The study of such examples will be covered in a subsequent paper by the authors.  Many of the results on the triadic Cantor set hold for a much larger class of examples.  In particular, it can be shown that every attractor of a self-similar iterated function system that satisfies the strong separability condition is such that its natural metric coming from the attractor is equivalent to a regular ultrametric.  This result also holds for cookie-cutter systems which is a class of Cantor sets that includes many Julia sets.  Basic definitions of these two classes of Cantor sets can be found in \cite{FAL2}.  An important generalization by the authors of some of the results for the triadic Cantor set is the following.

\begin{theo}
\label{cantor07.th-hmeas}
Let $C$ be the attractor of a self-similar iterated function system that satisfies the strong separability condition.  Then the following are true:

\noindent (i) $C$ is a $\zeta$-regular Cantor set with respect to a regular ultrametric that is metrically equivalent to the natural metric coming from the iterated function system;

\noindent (ii) up to a constant, $\mu$ is equal to the $s_0$-Hausdorff measure where $s_0$ is the similarity dimension of $C$.
\end{theo}

\noindent It is unclear whether $\zeta$-regularity is enough of a constraint in general to guarantee that the Hausdorff dimension and upper box dimension coincide.  Moreover, it is an open problem to find conditions under which the measure $\mu$ of a $\zeta$-regular Cantor set is actually the Hausdorff measure.

\noindent Another important class of examples is given by the transversal $\Xi$ of an aperiodic, repetitive Delone set of finite type \cite{BHZ}. Such an example can be given a natural tree structure coming from its patches. The Voronoi metric is then a natural regular ultrametric on $\Xi$. The special case of the Fibonacci tiling has been investigated by the authors. It can be shown that it is a $\zeta$-regular Cantor set with $\zeta$-function equal to the Riemann $\zeta$-function plus a small perturbation.  Because the Riemann $\zeta$-function has an isolated pole at $z=1$, then the Fibonnaci tiling has upper box dimension equal to $1$.  The algorithmic complexity of the Fibonnaci tiling is also $1$ and it seems that the upper box dimension and algorithmic complexity should agree for more general tilings. In \cite{LAP}, Lapidus proposes a new definition of fractality as a set whose $\zeta$-function has nonreal singularities in the positive half-plane. The Fibonacci tiling then provides a counterexample to this definition since it has only a singularity at $z=1$ in the positive half-plane. The Fibonacci tiling also has a natural construction as a cut-and-project tiling. The transversal space of the cut-and-project tiling gives a natural embedding (but not an isometry) of the transversal of the Fibonacci tiling into $\RM$. The authors have then shown that the measure $\mu$ associated to the $\zeta$-function of $\Xi$ is then the Lebesgue measure coming from this embedding. For this reason one can argue whether the transversal of the Fibonacci tiling is really a fractal. A generalization of this fact to the transversal of a cut-and-project tiling is a subject of future research.

\vspace{.5 cm}

\appendix
\section{Proof of Michon's Correspondence}
\label{cantor07.sect-appendix}

\vspace{.3 cm}
\subsection{Proof of Proposition \ref{cantor07.prop-umpfs}}
\label{cantor07.ssect-pfumpfs}

\noindent  
Given a regular ultrametric $d$, the equivalence relation $\stackrel{\epsilon}{\sim}$ given by $\epsilon$-chains will be shown to be a profinite structure.  (i)For $y\in [x]_\epsilon$, $B_{\epsilon}(y):=\{z\in C: d(z,y) < \epsilon\} \subset [x]_\epsilon$.  Thus $[x]_\epsilon$ is open.  Therefore $R_\epsilon=\bigcup_{x\in C} [x]_\epsilon\times[x]_\epsilon$ is open.  A compact metric space is totally bounded, so there exists $\epsilon$ such that $R_\epsilon=C\times C$.  

\noindent(ii)Let $x\stackrel{\epsilon}{\sim} y$.  Then there exists $x_0=x,x_1,\dots,x_n=y$ with $d(x_i,x_{i+1})<\epsilon$.  If $\eta=(\max\{d(x_i,x_{i+1}):0\leq i < n\})/2$ then $x\stackrel{\eta}{\sim} y$ with $\eta < \epsilon$.

\noindent(iii)Suppose $[x]_0:=\bigcap_{\epsilon\in\RM^+} [x]_\epsilon$ is the disjoint union of two closed sets $U$ and $V$. Since $C$ is compact, if both $U$ and $V$ are nonempty then there exists $u\in U$ and $v\in V$ such that $\tdist(U,V)=d(u,v) > 0$.  But then if $\eta=d(u,v)/2$ then $u\stackrel{\eta}{\nsim} v$.  So $[x]_0$ must be connected.  Thus since $C$ is totally disconnected, $[x]_0=\{x\}$.  Therefore, $\bigcap_{\epsilon\in\RM^+} R_\epsilon = \Delta$.  

\noindent Finally, given another regular ultrametric $d'\neq d$ then there exists $x,y\in C$ with $d(x,y)\neq d'(x,y)$.  Suppose that $d(x,y)=\epsilon > d'(x,y)=\epsilon'$.  If $\eta = (\epsilon+\epsilon')/2$, then $x \stackrel{\eta}{\sim_{d'}} y$ but $x \stackrel{\eta}{\nsim_{d}} y$ and therefore they give different profinite structures.

\noindent Conversely, given a profinite structure $\{R_\epsilon:\epsilon\in\RM^+\}$ on $C$ let $d(x,y):=\inf\{\epsilon: x\stackrel{\epsilon}{\sim} y\}$.  That $d(x,y)=0$ if and only if $x=y$ follows from the fact that $\bigcap_{\epsilon\in\RM^+} R_\epsilon = \Delta$.  For $x,y,z\in C$, if $x \stackrel{\epsilon_1}{\sim} y$ and $y \stackrel{\epsilon_2}{\sim} z$ and if $\epsilon=\max\{\epsilon_1,\epsilon_2\}$, then $x \stackrel{\epsilon}{\sim} z$.  Thus $d(x,z)\leq \max\{d(x,y),d(y,z)\}$ and $d$ is an ultrametric.  In order to show that $d$ is regular, let $id:C\to C$ be the identity map from $C$ with the original topology to $C$ with the metric topology.  First of all, if $x \stackrel{a}{\sim} y$ then by (ii) $x \stackrel{a-\delta}{\sim} y$ for some $\delta > 0$ and $d(x,y) < a$.  Thus, $d(x,y) < \epsilon$ if and only if $x \stackrel{\epsilon}{\sim} y$.  This gives that $B_a(x)=[x]_a$.  In fact, $[x]_\epsilon$ is open in the original topology.  This can be seen as follows.  Let $(x,y)\in C\times C$.  Since $R_\epsilon$ is open, then there exists an open set $V\subset C\times C$ such that $(x,y)\in V\subset R_\epsilon$.  But $C\times C$ has the product topology and therefore there exists open sets $U_x,U_y\subset C$ such that $(x,y)\in U_x\times U_y\subset V$.  For any $y\in U_y$, $(x,y)\in R_\epsilon$ and consequently $U_y\subset [x]_\epsilon$ and $[x]_\epsilon$ is open.  Therefore, $id$ is a continuous, bijective map from a compact space to a Hausdorff space and therefore a homeomorphism.  Thus, $d$ is regular.

\noindent Given two different profinite structures $\{R_\epsilon\}$ and $\{R'_\epsilon\}$, then without loss of generality there exists $\epsilon > 0$ and $(x,y)\in R_\epsilon$ such that $(x,y)\notin R'_\epsilon$.  Suppose $\{R_\epsilon\}$ gives ultrametric $d$ and $\{R'_\epsilon\}$ gives ultrametric $d'$.  Then by (ii), $(x,y)\in R_{\epsilon-\delta}$ for some $\delta > 0$ and $d(x,y) < \epsilon \leq d'(x,y)$.  Consequently, $d\neq d'$.
\hfill $\Box$

\vspace{.3 cm}
\subsection{Proof of Proposition \ref{cantor07.prop-wrtpfs}}
\label{cantor07.ssect-pfwrtpfs}

\noindent  
Let $d$ be a regular ultrametric on $C$ and let $\{R_\epsilon \}$ be the profinite structure corresponding to $d$.  The tree $\Tt$ is built as follows.  Let $\epsilon_0=\inf\{\epsilon:R_\epsilon=C\times C\}$.    Then $R_{\epsilon_0}\neq C\times C$ since $R_{\epsilon_0}=\bigcup_{\epsilon' < \epsilon_0}R_\epsilon$.  Similarly, let $\epsilon_{i+1}=\inf\{\epsilon:R_\epsilon=R_{\epsilon_i}\}$.  Then $\{\epsilon_i\}_{i=0}^\infty$ is such that $R_{\epsilon_i}\neq R_{\epsilon_{i+1}}$.  Let the root of $\Tt$ correspond to $C$ and let the vertices of height $n$ correspond to the equivalence classes of $R_{\epsilon_{n-1}}$.  Let the edges be defined by $[x]_{\epsilon_j} \succeq [y]_{\epsilon_k}$ if and only if $[x]_{\epsilon_j} \supset [y]_{\epsilon_k}$.  Then $\Tt$ is a rooted tree with no dangling vertex.  As seen in the proof of the previous proposition, every equivalence class is clopen.  Thus each vertex has a finite number of children and has a descendant with more than one children.  So, $\Tt$ is a Cantorian tree.  In general, $\Tt$ is not reduced.  However, since each vertex has a descendant with more than one child, edge reduction can be applied to each vertex with only one child without altering $\pT$.  This will give a reduced tree $\Tt'$ with vertices $\Vv'\subset \Vv$ such that $\pT'=\pT$ as topological spaces.

\noindent Let $\Phi:\pT'\to C$ be defined by $\Phi(v_0v_1\cdots)=\bigcap_{i=1}^\infty [x_i]_{\epsilon_i}$ where $v_i=[x_i]_{\epsilon_i}$.  This map is bijective and $\Phi^{-1}([x]_{\epsilon_i})=[v]$ where $v=[x]_{\epsilon_i}$.  Thus $\Phi$ is continuous and since $\pT'$ is compact, $\Phi$ is a homeomorphism.  By abuse of notation, let $[v]=[x]_{\epsilon_i}$if $v=[x]_{\epsilon_i}$. 

\noindent If $v=[x]_{\epsilon_k}$ then let $\epsilon(v):=\epsilon_{k+1}$.  Since $\epsilon_k > 0$ for all $k$, then $\epsilon:\Vv'\to\RR^+$.  (i) follows automatically.  (ii)  Since $\epsilon([x]_{\epsilon_k})\leq \epsilon_k$ and $\epsilon_k\to 0$ then $\lim_{k\uparrow\infty} \epsilon([x]_{\epsilon_k})\leq \lim_{k\to\infty} \epsilon_k=0$.  So $\Tt'$ is a reduced, weighted, rooted Cantorian tree.

\noindent Let $\Tt$ be a reduced, rooted Cantorian tree with weight function $\epsilon$.  For $x,y\in \pT =: C$, let $d(x,y)=\epsilon(x\wedge y)$ for $x\neq y$ and $d(x,x)=0$.  It is straightforward to show that $d$ is an ultrametric on $C$.  Given $r > 0$ and $x\in C$, let $B_r(x):=\{y\in C: d(x,y) < r\}$.  By (ii), $B_r(x)$ has more than one point, so let $v=l.c.p.(B_r(x))$.  By the definition of $v$, for $y\in [v]$ there exists $z\in B_r(x)$ such that $x\wedge y\preceq x\wedge z$.  Thus $d(x,y)\leq d(x,z) < r$ and therefore $[v]=B_r(x)$.  Consequently, $B_r(x)$ is open in $\pT$ and $d$ is regular.  

\noindent For $x,y\in [v]$ then $x\wedge y \preceq v$ and $d(x,y)=\epsilon(x\wedge y)\leq \epsilon(v)$.  Thus, $\tdiam([v])\leq \epsilon(v)$.  Conversely, since $v$ has more than one child then there exists $x,y\in [v]$ such that $v=x\wedge y$.  Therefore, $\epsilon(v)=d(x,y)\leq \tdiam(v)$ and $\epsilon(v)=\tdiam([v])$.

\noindent Starting with a regular ultrametric $d$ on $C$, let $d_\epsilon$ be the regular ultrametric obtained from the Cantorian tree $\Tt$ corresponding to $d$.  Let $x,y\in C$.  Then $d_\epsilon(x,y)=\epsilon(x\wedge y)=\epsilon_{k+1}$ if $x\wedge y=[x]_{\epsilon_{k}}$.  So $x\stackrel{\epsilon_{k+1}}{\nsim} y$ but $x\stackrel{\epsilon_{k+1}+\delta}{\sim} y$ for $\delta > 0$.  Since $d$ is an ultrametric then $d(x,y)=\epsilon_{k+1}$.  Thus $d=d_\epsilon$ and $\pT$ is isometric to $C$.

\noindent Starting with a reduced, weighted, rooted Cantorian tree $\Tt$ let $\Tt_d$ be the tree obtained from the regular ultrametric $d$ corresponding to $\Tt$.  Let $\Phi$ be the homeomorphism from $\pT\to\pT_d$.   Let $\Psi:\Vv\to\Vv_d$ be defined by $\Psi(v)=l.c.p(\Phi([v]))$.  Because each tree is reduced there is a one-to-one correspondence between clopen sets in the boundary and vertices, thus $\Psi$ is a bijection.  Therefore the correspondence between reduced, weighted, rooted Cantorian trees and regular ultrametrics is indeed a bijection.
\hfill$\Box$

\vspace{.3cm}

\newpage


\begin{thebibliography}{BEL}

\bibitem {AlK} S. Albeverio and W. Karwowski, {A Randow Walk on $p$-adics - the generator and it spectrum}.  {\sl Stoch. Process. Appl.} {\bf 53} (1994), 1-22.

\bibitem{BBG} J. Bellissard; R. Benedetti; J.-M. Gambaudo, { Spaces of Tilings, Finite Telescopic Approximations and Gap-labelling}, {\sl Commun. Math. Phys.} {\bf 261} (2006), 1-41.

\bibitem{Bel} J. Bellissard,  {Noncommutative Geometry of Aperiodic Solids}. In {\sl Geometric and Topological Methods for Quantum Field Theory: Proceedings of the Summer School Held in Villa de Leyva, Colombia, July 9-27 2001}, World Scientific, River Edge, NJ, 2003, 86-156.

\bibitem{BHZ} J. Bellissard; D. Hermmann; M. Zarrouati, Hull of Aperiodic Solids and Gap Labeling Theorems. In {\sl Directions in Mathematical Quasicrystals}, CRM Monograph Series 13, AMS, Providence, RI, 2000, 207-259.

\bibitem{Bla86} B. Blackadar,  {\sl $K$-Theory for Operator Algebras}. Cambridge University Press 1998.

\bibitem{BOL} B. Bollobas, {\sl Modern Graph Theory}. Springer-Verlag 1998.

\bibitem{ChI} E. Christensen and C. Ivan, {Sums of Two Dimensional Spectral Triples}. {\sl Math. Scand.} (1) {\bf 100} (2007),   35-60.

\bibitem{CIL} E. Christensen; C. Ivan; M. Lapidus,  {Dirac Operators and Spectral Triples for some Fractal Sets Built on Curves}.  {\tt arXiv:math/0610222 v2 (2007)}.

\bibitem{CO88} A. Connes, {Trace de Dixmier, modules de Fredholm et G\'eom\'etrie riemannienne}. (French) [Dixmier trace, Fredholm modules and Riemannian geometry]. In {\sl Conformal field theories and related topics},  (Annecy-le-Vieux, 1988),  Nuclear Phys. B Proc. Suppl.  5B, 1988,   65-70.

\bibitem{CON} A. Connes, {\sl Noncommutative Geometry}. Academic Press 1994.

\bibitem{CONW} J. Conway, {\sl A Course in Functional Analysis}. Springer-Verlag 1990.

\bibitem{dMF} M. Del Muto and A. Fig\`a-Talamanca, {Diffusion on Locally Compact Ultrametric Spaces}. {\sl Expo. Math.} (3) {\bf 22} (2004), 197–211.

\bibitem {Eva} S. Evans,  {Local fields, Gaussian measures, and Brownian motions}.  In {\sl Topics in probability and Lie groups: boundary theory}.  CRM Proc. Lecture Notes 28, Amer. Math. Soc., Providence, RI, 2001, 11-50.

\bibitem {FaJ} C. Favre and M. Jonsson  {\sl The Valuative Tree}.  Springer-Verlag 2004.

\bibitem{FAL} K. Falconer, {\sl Fractal Geometry: Mathematical Foundations and Applications}. John Wiley and Sons 1990.

\bibitem{FAL2} K. Falconer, {\sl Techniques in Fracal Geometry}.  John Wiley and Sons 1997.

\bibitem{Fig} A. Fig\`a-Talamanca,  {Diffusion on compact ultrametric spaces}. In {\sl Noncompact Lie groups and some of their applications} (San Antonio, TX, 1993), NATO Adv. Sci. Inst. Ser.C Math. Phys. Sci. 429, Kluwer Acad. Publ., Dordrecht, 1994, 157–167.

\bibitem{Fuk} M. Fukushima,  {\sl Dirichlet Forms and Markov Processes}.  North-Holland 1980.

\bibitem{GRA} J. Gracia-Bond\'ia; J. V\'arilly; H. Figueroa, {\sl Elements of Noncommutative Geometry}.  Birkhauser 2001.

\bibitem{GuI1} D. Guido and T. Isola,  {Dimensions and Singular Traces for Spectral Triples, with Applications to Fractals}. {\sl J. Funct. Anal.} {\bf 203} (2003), 362-400.

\bibitem{GuI2} D. Guido and T. Isola,  {Dimension and Spectral Triples for Fractals in $\RR^N$}.  In {\sl Advances in Operator Algebras and Mathematical Physics}, Theta Ser. Adv. Math. 5, Theta, Bucharest, 2005,   89-108.

\bibitem{HAU} F. Hausdorff, {\sl Set Theory} (Translated in English by John R. Aumann). Chelsea Publishing Co. 1962.

\bibitem{HaR} G.H. Hardy and M. Riesz, {\sl The General Theory of Dirichlet's Series}.  Cambridge University Press 1915.

\bibitem{KIG} J. Kigami, {\sl Analysis on Fractals}.  Cambridge University Press 2001.

\bibitem{Koc} A. Kochubei,  {\sl Pseudo-Differential Equations and Stochastics Over Non-Archimedan Fields}.  Marcel Dekker Inc. 2001.

\bibitem{Koz} S.V. Kozyrev, {Wavelet Analysis as a $p$-adic Spectral Analysis}. {\sl Izv. Ross. Akad. Nauk Ser. Mat.} {\bf 66} (2002),   149-158.

\bibitem{LAP} Lapidus, Michel. {Towards a Noncommutative Fractal Geometry? Laplacians and Volume Measures on Fractals}. In {\sl Harmonic Analysis and Nonlinear Differential Equations}, Contemp. Math. 208, Amer. Math. Soc., Providence, RI, 1997, 211-252.

\bibitem{MIC} G. Michon, { Les Cantors r\'eguliers}. {\sl C. R. Acad. Sci. Paris S\'er. I Math.} (19) {\bf 300} (1985), 673--675.

\bibitem{Par} K.R. Parthasarathy,  {\sl Probability Measures on Metric Spaces}.  Academic Press 1967.

\bibitem{Rud} W. Rudin,  {\sl Real and Complex Analysis}.  McGraw-Hill 1987.

\bibitem{Sch} W.H. Schikhof,  {\sl Ultrametric Calculus: An introduction to $p$-adic analysis}.  Cambridge University Press 1984.

\bibitem{VVZ} V.S. Vladimirov; I.V. Volovich; E.I. Zelenov, {\sl $p$-adic Analysis and Mathematical Physics}.  World Scientific 1994.

\vspace{1cm}


\end{thebibliography}
\end{document}